\documentclass[11pt]{article} 
\usepackage{amsmath,amsthm,amsfonts,latexsym,amssymb,enumerate,color}
\usepackage{graphicx}
\usepackage[font=small,labelfont=bf]{caption}

\def\CC{\mathbb C}
\def\CP{{\rm CP}}
\def\D{\mathcal D}
\def\F{\mathcal F}
\def\FBP{{\rm FBP}}
\def\DD{\mathbb D}
\def\G{\mathcal G}
\def\HH{\mathcal H} 
\def\K{\mathcal K}

\def\N{\mathcal N}
\def\R{\mathcal R}
\def\RR{\mathbb R}
\def\SS{\mathcal S} 
\def\TT{\mathbb T}
\def\ZZ{\mathbb Z}

\def\phi{\varphi}

\def\T{Toeplitz } 
\def\with{\qquad \hbox{with} \quad} 

\def\H2p{\overline{H^2_0}}

\def\beq{\begin{equation}}
\def\eeq{\end{equation}}
\def\ds{\displaystyle}

\def\diag{\mathop{\rm diag}\nolimits}
\def\gcd{\mathop{\rm gcd}\nolimits}

\def\ind{\mathop{\rm ind}\nolimits}
\def\ker{\mathop{\rm ker}\nolimits}

\def\spam{\mathop{\rm span}\nolimits} 

\newtheorem{thm}{Theorem}[section]
\newtheorem{prop}[thm]{Proposition}

\newtheorem{lem}[thm]{Lemma}
\newtheorem{cor}[thm]{Corollary}
\newtheorem{rem}[thm]{Remark}

\newtheorem{ex}[thm]{Example}

\numberwithin{equation}{section}

\def\beginpf{\begin{proof}}
\def\endpf{\end{proof}}

\begin{document}

\title{Scalar-type kernels for block Toeplitz operators}

\date {}

\author{M.~Cristina C\^amara\thanks{
Center for Mathematical Analysis, Geometry and Dynamical Systems,
Instituto Superior T\'ecnico, Universidade de Lisboa, 
Av. Rovisco Pais, 1049-001 Lisboa, Portugal.
 \tt ccamara@math.ist.utl.pt}\ \thanks{Corresponding author} \ and  Jonathan R.~Partington\thanks{School of Mathematics, University of Leeds, Leeds LS2~9JT, U.K. {\tt j.r.partington@leeds.ac.uk}}
}

\maketitle

\begin{abstract}
It is shown that the kernel of a Toeplitz operator with $2\times 2$ symbol $G$ can be described exactly in terms of any given function in a very  wide class, its image under multiplication by $G$, and their left inverses, if the latter exist. As a consequence, under many circumstances the kernel of a block Toeplitz operator may be described as the product of a
space of scalar complex-valued functions by a fixed column vector of functions. Such kernels are said to 
be of scalar type, and in this paper they are studied and described explicitly in many concrete situations. Applications are given to
the determination of kernels of truncated Toeplitz operators
for several new classes of symbols.
\end{abstract}



\maketitle

\section{Introduction}

Kernels of Toeplitz operators (also called Toeplitz kernels) have generated an enormous interest for various reasons, among which is the fact
that they have fascinating properties and a rich structure, they are important in many applications,
and several relevant classes of analytic functions can be presented as kernels of Toeplitz operators.
For instance, model spaces (defined below) are Toeplitz kernels.
Two recent surveys of this area are \cite{HM} and \cite{CPiwota}.

It is natural to expect that kernels of block Toeplitz operators, whose study provides a clear example of the fruitful interplay between operator theory, complex analysis and linear algebra, will have an even richer and more involved structure. The case of Toeplitz operators with $2\times 2$ symbols is particularly interesting for its connections to truncated Toeplitz operators \cite{CP17,CP16} and the corona theorem \cite{CDR},
and because their study leads to various surprising results.
One of these unexpected results is that, as we prove in Theorem \ref{thma:3.1}, one can explicitly describe the kernel of a Toeplitz operator with $2\times 2$ symbol $G$, not necessarily bounded, in terms of any given function $f$ in a very wide class, its image under multiplication by $G$, and their left inverses, 
assuming that the latter exist. As a result of this, we show that, although those kernels consist of vector functions, in many cases they behave as having a scalar nature, since they can
be expressed as the product of a space of scalar functions by a fixed vector function. These
kernels will be called {\em scalar-type Toeplitz kernels}.

A natural question arising in this case regards which properties of scalar Toeplitz kernels remain valid for scalar-type Toeplitz kernels. While Coburn's Lemma, stating that for $\phi\in L_\infty(\mathbb R)$, either $\ker T_\phi$ or $\ker T_\phi^*$ is zero, cannot be extended in the same form to the case of a general $2\times  2$ bounded symbol $G$, even if  $\ker T_G$ and $\ker T_G^*$ are both of scalar type,  Theorem \ref{theor3.7} may be seen as a version of Coburn's Lemma for $2\times 2$ symbols. Moreover we show that any scalar-type Toeplitz kernel is the product of a fixed vector function by a scalar nearly $S^*$-invariant space $\K$, which is closed if $G\in L_\infty^{2\times 2}$ and therefore, using a well known result by Hitt \cite{Hitt}, can be characterized as the product of a scalar Toeplitz kernel, in fact a model space, by a fixed $2\times 1$ function (Theorem \ref{thma:3.15}). Although this model space is not known in general, by using the corona theorem we obtain sufficient conditions for $\K$ to be a model space, explicitly described in terms of the functions $f$ and $g$, leading to conditions for injectivity and invertibility for $T_G$ (Theorem \ref{thma:3.12} and its corollaries). 
We note that some related results can be found in \cite[Prop. 4.6]{FHR2}.
We show moreover that, as in the case of scalar symbols, every scalar-type Toeplitz kernel has a maximal function (Theorem \ref{thma:3.17}).

The results of Section 3 are applied in Section 4 to study and describe the kernel of truncated Toeplitz operators in two different classes which extend previously studied ones. In the first case we show that the kernels are given by the product of a model space, which is explicitly determined, by a fixed vector function in $H_\infty^+$, and we establish necessary and sufficient conditions for injectivity and invertibility of the truncated Toeplitz operators. In the second case we also obtain an explicit description of the kernel as a product of a scalar Toeplitz kernel by a fixed function.

\vspace{5mm}

We write $\CC^+$ and $\CC^-$ for the upper and lower complex half-planes, and
$H_p^\pm$ $(1 \le p \le \infty)$ for the associated Hardy spaces of analytic functions on $\CC^\pm$.
The operators $P^\pm$ are the standard Riesz projections from $L_p=L_p(\RR)$ onto
the subspaces $H_p^\pm$.
It will be recalled that functions in $H_p^+$ have inner/outer factorizations. For two inner
functions $\theta, \phi$ we write $\theta \preceq \phi$ or $\phi \succeq \theta$
to mean that $\theta$ is a divisor of $\phi$
in $H_\infty^+$. We may also use the strict versions of these relations, written $\prec$ and $\succ$.

The Smirnov class $\N_+$ consists of all analytic functions $f=g_+/h_+$, where 
$g_+ \in H_1^+$ and $h_+ \in H_2^+$  
with $h_+$ outer. We may instead
take $g_+, h_+ \in H_\infty^+$ (see, e.g. \cite{nik}).

These notions can be found in standard texts on Hardy spaces, such as \cite{koosis} and \cite{nik}.

For $G \in L_\infty^{n \times n}$, with $n=1,2,\ldots$, the Toeplitz operator $T_G$ on $(H_2^+)^n$
is the composition $P^+ M_G$, where $M_G$ denotes multiplication by $G$.
For $\theta \in H_\infty^+$ inner, the model space $K_\theta$ is $\ker T_{\overline\theta}$, which
equals $H_2^+ \ominus \theta H_2^+=H_2^+ \cap \theta H_2^-$.

For a unital algebra $\mathcal A$, we write $\G({\mathcal A})$ for the group of invertible elements.

We use the notation $(f,g)$ interchangeably with $[f \; g]^T=\ds \begin{bmatrix}f \\ g \end{bmatrix}$.


\section{Motivation: matrix symbols with a bounded factorization}
\label{sec:2}

Let $G \in (L^\infty)^{2 \times 2}$ admit a bounded (Wiener--Hopf) factorization 
(\cite{CG,MP}) on the real line, of the form
\beq\label{eq:2.1} 
G=G_- \diag(r^{k_1},r^{k_2}) \, G_+^{-1},
\eeq
where $G_\pm \in \G(H_\infty^\pm)^{2 \times 2}$, $k_1, k_2 \in \ZZ$, and
\beq\label{eq:2.2}
r(\xi)= \frac{\xi-i}{\xi+i}, \qquad \hbox{for} \quad \xi \in \RR.
\eeq
The class of matrix functions admitting such a factorization includes, in particular, all $2 \times 2$ matrix functions $G$ with elements in the algebra $C^\mu(\dot \RR)$ of H\"older-continuous functions
in $\dot\RR:= \RR \cup \{\infty\}$ with exponent $\mu \in (0,1)$, 
 or in the Wiener algebra $W(\dot \RR)$
(\cite{CG,MP}),
as long as $\det G \in \G L^\infty$. If $k_1=k_2 \geq 0$ in \eqref{eq:2.1}, then $\ker T_G=\{0\}$;
if at least one of the integers $k_1,k_2$ is negative, then
$\ker T_G \ne \{0\}$. Let us assume, for simplicity, that
$\ind(\det G)=0$, in which case $k_1=-k_2=-k$, say, and \eqref{eq:2.1} takes the form
\beq\label{eq:2.3}
G=G_- \diag(r^{-k},r^{k}) G_+^{-1}, \qquad k \in \ZZ_0^+=\{0,1,2,\ldots\},
\eeq
and let $G_\pm = \left[ g^\pm_{ij} \right]_{i,j=1,2}$. Rewriting the equation \eqref{eq:2.3} as
\beq\label{eq:2.4}
GG_+\diag(r^k, r^{-k}) = G_-
\eeq
and taking the first columns of the matrices on the left and right-hand side of \eqref{eq:2.4}, we obtain
\beq\label{eq:2.5}
Gr^k g_+=g_-, \qquad \hbox{with}\quad  g_\pm=\left( g^\pm_{11}, g^\pm_{21}\right).
\eeq
On the other hand, $\ker T_G$ consists of all functions $\phi_+ \in (H_2^+)^2$ such that
\beq\label{eq:2.6}
G\phi_+=\phi_- \qquad \hbox{with} \quad \phi_- \in (H_2^-)^2.
\eeq
From \eqref{eq:2.5} and \eqref{eq:2.6} we have then
\beq\label{eq:2.7}
G[r^k g_+ \quad \phi_+]=[g_- \quad \phi_-]
\eeq
and, on taking determinants on both sides and noting that $\det G=d_-d_+^{-1}$ where $d_\pm = \det G_\pm \in \G H_\infty^\pm$,
it follows that 
\beq\label{eq:2.8}
d_+^{-1}\det[r^k g_+ \quad \phi_+]=d_-^{-1}\det[g_- \quad \phi_-].
\eeq
The left-hand side of this identity is in $H_2^+$, while the right-hand side is in $H_2^-$.
Consequently, they are both equal to zero and we have that
\beq\label{eq:2.9}
\phi_+=\Lambda r^k g_+, \qquad \phi_-= \tilde\Lambda g_-,
\eeq
where $\Lambda$ and $\tilde\Lambda$ are scalar functions defined a.e.\ on $\RR$.
To show that $\Lambda = \tilde\Lambda$, we now take into account the fact that
the invertibility of $G_\pm$ in $(H_\infty)^{2 \times 2}$ means that
the first column of $G_\pm$ is a corona pair in $\CC^\pm$ (\cite{Sarason,CDR}) and therefore left-invertible in $H_\infty^\pm$, 
with left inverse given by $\tilde g^T_\pm$, where
$\tilde g_\pm = (g_{22}^\pm, -g^\pm_{12})/d_\pm$,
and
\beq\label{eq:2.10}
\tilde g^T_\pm g_\pm = 1 \qquad \hbox{in} \quad \CC^\pm.
\eeq
In fact, from \eqref{eq:2.6} and \eqref{eq:2.9} we have
\[
G(\Lambda r^k g_+)=\tilde\Lambda g_-,
\]
and it follows from \eqref{eq:2.5}
that $\Lambda g_- =\tilde \Lambda g_-$. Multiplying both sides by $\tilde g^T_-$ on the left we conclude that
$\Lambda = \tilde\Lambda$. Then multiplying both sides of the two equations in \eqref{eq:2.9}
by $\tilde g^T_+$ and $\tilde g^T_-$ respectively, we obtain moreover that
\beq\label{eq:2.11}
\Lambda = r^{-k}(\tilde g_+^T\phi_+)=\tilde g^T_-\phi_-  \with \tilde g_\pm^T\phi_\pm \in H_2^\pm.
\eeq
Therefore $r^k\Lambda = \tilde g_+^T\phi_+ \in \ker T_{r^{-k}}=K_{r^k}$,
where $K_{r^k}$ is the model space $H^+_2 \ominus r^k H_2^+$, and it follows from \eqref{eq:2.9} and \eqref{eq:2.11} that
\[
\ker T_G \subset K_{r^k}g_+.
\]
Conversely, $K_{r^k}g_+ \subset \ker T_G$ because if $P \in K_{r^k}$ then
\begin{eqnarray*}
G(Pg_+) &=& G_- \diag(r^{-k},r^k) G_+^{-1}(Pg_+)\\
&=& G_- \diag(r^{-k}P, r^k P)(1,0) \\
&=& G_-(r^{-k}P,0) \in (H_2^-)^2.
\end{eqnarray*}
Thus
\beq\label{eq:2.12}
\ker T_G = \K g_+,
\eeq
where $\K$ is a scalar model space, associated with the inner function $r^k$, and $g_+$ is a fixed vector function.
So we see that for a wide class of Toeplitz operators with $2 \times 2$ matrix symbols, including for instance all invertible
$2 \times 2$ H\"older-continuous matrices $G$ with $\ind(\det G)=0$, the
corresponding kernels are spaces of vector functions which can
nonetheless be described as the product of a certain space $\K$ of scalar functions by a fixed vector function.
We say in this case that $\ker T_G$ is a {\em scalar-type Toeplitz kernel}.

The same result  would hold in the case of any $2\times 2$ matrix symbol $G$ for which one can find a solution to
\beq\label{eq:2.13}
Gf=g
\eeq
with $f=r^kg_+$ and $g=g_-$ where $k \in \ZZ^+_0$, $g_\pm \in (H_\infty^\pm)^2$,
such that $g_\pm$ satisfy the condition of Carleson's corona theorem in $\mathbb C^\pm$
(\cite{garnett,CDR}).
These conditions can be seen in terms of left invertibility of $f$ and $g$ by saying that there
exist vector functions $\tilde g_\pm \in (H_\infty)^2$ such that
\beq\label{eq:2.14}
\tilde g^T_\pm g_\pm = 1 \qquad \hbox{in}\quad \CC^\pm.
\eeq
The main difficulty in applying these results consists in the fact that it is in general very difficult, or even impossible, to find solutions 
to \eqref{eq:2.13} satisfying the above-mentioned conditions. We may however find other solutions to
\eqref{eq:2.13},
satisfying less restrictive conditions; it is natural to ask then whether such a relation would still allow
us to describe the kernel of $T_G$, and whether the kernel would be of scalar type.

In the next section we shall show that it is indeed possible to describe the kernel of a Toeplitz operator with $2 \times 2$ symbol
$G$ and give conditions for it to be a scalar-type kernel, for a very general set of symbols, in terms of a solution 
to $Gf=g$ where $f$ and $g$ are assumed to be left-invertible vector functions in a very general class.
In particular, we shall not assume any analyticity conditions on the functions $f$ and $g$.

\section{Scalar-type kernels for \T operators with $2 \times 2$ symbols}
\label{sec:3}

Let $\F$ denote the   space of all complex-valued functions defined almost everywhere on $\RR$,
where as usual we identify two functions if they are equal almost everywhere. Let $G \in \F^{2 \times 2}$
and let
\beq
\D= \{ f_+ \in (H_2^+)^2: G f_+ \in (L_2)^2 \}.
\eeq
The operator $T_G: \D \to (H_2^+)^2$ defined by
\beq
T_G \,f_+ = P^+ (Gf_+), \quad f_+ \in \D,
\eeq
where $P^+: (L_2)^2 \to (H_2^+)^2$ denotes the orthogonal projection, 
is called the {\em Toeplitz operator} with symbol $G$. 
If $G \in (L_\infty)^{2\times 2}$, then $T_G$ is a 
bounded operator on $(H_2^+)^2$. 

Another class of symbols of interest arises on taking
$G \in \lambda_+ (L_2)^{2\times 2}$, where $\lambda_+(\xi)=\xi+i$; then
$T_G$ is densely defined on $(H_2^+)^2$. These symbols correspond to $(L^2(\TT))^{2\times 2}$
symbols of Toeplitz operators on
$(H^2(\DD))^2$    (here $\DD$ is the unit disc and $\TT$  the unit circle in $\CC$).

In what follows $f$ and $g$ denote left-invertible functions in $\F^{2 \times 1}$ with left inverses
$\tilde f^T$ and $\tilde g^T$, where $\tilde f, \tilde g \in \F^{2 \times 1}$. We assume moreover that $G \in \F^{2 \times 2}$ and, unless said otherwise, $\det G \in \G \F$.

We write $J$ for the matrix $\ds J=\begin{bmatrix}0 & -1 \\ 1 & 0 \end{bmatrix}$.

The following result shows that, surprisingly, given a Toeplitz operator $T_G$, one can describe its kernel in terms of any given function $f$, its image under multiplication by $G$, $g=Gf$, and their left inverses if they exist. Although, necessarily, with a somewhat technical appearance, it leads on to explicit characterizations of Toeplitz kernels, as we shall show here.

\begin{thm}\label{thma:3.1}
If $f$ and $g$ are left-invertible functions in $\F^{2 \times 1}$ such that
\beq\label{eqa:3.1}
Gf=g
\eeq
and we define 
\beq\label{eqa:3.2}
\SS = \left( \det G . f^T (H_2^+)^2 \right) \cap \left( g^T (H_2^-)^2 \right),
\eeq
and, for each $s \in \SS$,
\beq
\HH_+^s = \left\{ \psi_+ \in (H_2^+)^2: f \tilde f^T \psi_+ + \frac{s}{\det G} J\tilde f \in (H_2^+)^2 \right \}
\eeq
and
\beq
\HH^s_- =  \left\{ \psi_- \in (H_2^-)^2: g \tilde g^T \psi_- + {s} J\tilde g \in (H_2^-)^2 \right \},
\eeq
then $\phi_+ \in \ker T_G$ if and only if
\beq\label{eqb:3.5}
\phi_+ = \Lambda f + \frac{s}{\det G} J \tilde f,
\eeq
where 
\beq\label{eqa:3.6}
s \in \SS \hbox{ with } \HH^s_\pm \ne \emptyset, \quad \Lambda \in \tilde f ^T\HH^s_+, \quad \Lambda 
+ \frac{s}{\det G} \tilde g^T (GJ \tilde f) \in \tilde g^T \HH^s_-.
\eeq
\end{thm}

As we shall see later, it is possible for the space $\SS$ in the statement of the theorem to reduce to $\{0\}$,
in which case $\HH^s_\pm$ may be empty.

Naturally, although $f$ and $g$ in \eqref{eqa:3.1} can be chosen in a very general class, it is important that they are such that the description given by Theorem~\ref{thma:3.1} for the kernel of $T_G$ is useful, in the sense that it allows for a good understanding of the kernel.
There is no general method to obtain good solutions, in that sense,
and the choice of $f$ and $g$ must be made on a case by case basis. However, as we show in 
Example \ref{ex:3.3r} below, the great degree of freedom that we are allowed in that choice can permit us to
obtain a clear description of $\ker T_G$ in terms of
natural solutions to $Gf=g$ with no particular analytic properties.

\beginpf
(i) Let $\phi_+ \in \ker T_G$, i.e., $\phi_+\in (H_2^+)^2$ and $G\phi_+=\phi_- \in (H_2^-)^2$. From
this identity and \eqref{eqa:3.1} we have that
\beq
G [f \;\; \phi_+] = [g \;\; \phi_-].
\eeq
Let
\beq
s=\det G .\det [f \;\; \phi_+]= \det[g \;\; \phi_-].
\eeq
Now, since $\det[f \;\; J\tilde f]=1$, we have
\[
\det \left[f \; \;\frac{s}{\det G} J\tilde f \right]=\frac{s}{\det G}=\det[f \; \;\phi_+]
\]
and analogously
\[
\det [g \; \;sJ \tilde g] =s=\det[g \;\; \phi_-].
\]
Therefore,
\[
\det \left[f \quad \left( \phi_+-\frac{s}{\det G}J\tilde f \right)\right] = 0 = \det[g \quad (\phi_--sJ\tilde g)],
\]
and, since $f$ and $g$ are left invertible, it follows that there are scalar functions $\Lambda, \tilde \Lambda \in \F$ such that
\begin{eqnarray}
\phi_+ &=& \Lambda f + \frac{s}{\det G} J \tilde f, \label{eqa:3.9}
\\
\noalign{ \hbox{and}} \nonumber \\
\phi_- &=& \tilde \Lambda g + sJ \tilde g. \label{eqa:3.10}
\end{eqnarray}
Multiplying \eqref{eqa:3.9} and \eqref{eqa:3.10} on the left by $\tilde f^T$ and $\tilde g^T$, respectively, and taking
into account the fact that
$\tilde f^T J \tilde f = \tilde g^T J \tilde g=0$, we get
\beq\label{eqb:3.13}
\Lambda=\tilde f^T \phi_+, \qquad \tilde \Lambda = \tilde g^T \phi_-.
\eeq
Moreover, from \eqref{eqa:3.9}, \eqref{eqa:3.10}, and the assumption that $Gf=g$, we have
\beq
G \phi_+=\phi_- 
\implies \Lambda g + \frac{s}{\det G} GJ\tilde f = \tilde \Lambda g + sJ \tilde g,
\eeq
and, multiplying the last equation on the left by $\tilde g^T$, we get
\beq\label{eqa:3.12}
\Lambda + \frac{s}{\det G} (\tilde g^T GJ \tilde f)= \tilde\Lambda.
\eeq
On the other hand, multiplying \eqref{eqa:3.9} and \eqref{eqa:3.10} on the left by $f^T J$ and $g^TJ$ respectively, we have
\begin{eqnarray}
f^TJ \phi_+ &=& \frac{s}{\det G} f^T JJ\tilde f = - \frac{s}{\det G} f^T \tilde f = -\frac{s}{\det G}\,,
\\
\noalign{\hbox{and}} \nonumber \\
g^T J \phi_- &=& sg^T JJ \tilde g = -s,
\end{eqnarray}
taking into account the fact that $f^T \tilde f = g^T \tilde g = 1$. Therefore,
\[
s= -\det G . f^T J \phi_+ = - g^T J \phi_-,
\]
and we conclude that $s \in \SS$ and, since $\phi_\pm \in \HH^s_\pm$
by \eqref{eqa:3.9}--\eqref{eqb:3.13}, that
$\HH^s_\pm \ne \emptyset$.
From \eqref{eqa:3.9}--\eqref{eqa:3.12} we see that \eqref{eqa:3.6} holds.\\

(ii) Conversely, suppose that $\phi_+ = \Lambda f + \frac{s}{\det G} J \tilde f$, 
where $s, \Lambda$ and $\tilde \Lambda $ satisfy \eqref{eqa:3.6}. Then, since
$\Lambda \in \tilde f^T\HH^s_+$, we have, for some $\psi_+\in (H_2^+)^2$,
\[
\Lambda = \tilde f^T \psi_+, \quad \hbox{where } f\tilde f^T \psi_+ = - \frac{s}{\det G} J\tilde f + F_+, \quad F_+ \in (H_2^+)^2,
\]
and therefore
\[
\begin{array}{rclcl}
\phi_+ &=& \ds (\tilde f^T \psi_+) f + \frac{s}{\det G} J\tilde f 
= \ds f \tilde f^T \psi_+ + \frac{s}{\det G}J \tilde f \\
\\
&=& \ds - \frac{s}{\det G}J\tilde f + F_+ + \frac{s}{\det G}J\tilde f = \ds F_+ \in (H_2^+)^2.
\end{array}
\]
On the other hand, using Lemma \ref{lema:3.2} below, we have
\[
\begin{array}{rclcl}
G\phi_+& =& \ds \Lambda g +  \frac{s}{\det G} GJ\tilde f = \ds \Lambda g  + \frac{s}{\det G}[(\tilde g^T G J \tilde f) g + \det G. J \tilde g] \\
\\ 
&=& \ds \left(\Lambda + \frac{s}{\det G}(\tilde g^T GJ \tilde f)\right)g + sJ\tilde g  = (\tilde g^T \psi_-)g+ s J \tilde g= g\tilde g^T \psi_-+ s J \tilde g 
\end{array}
\]
with $\psi_- \in \HH^s_-$; therefore, $G \phi_+ \in (H_2^-)^2$ and it follows that $\phi_+ \in \ker T_G$.
\endpf

\begin{lem}\label{lema:3.2}
Let $Gf=g$; then $GJ\tilde f = \tilde g^T(GJ\tilde f)g + \det G. J\tilde g$.
\end{lem}
\beginpf
We have
$
G[f \;\; J\tilde f] = [g \;\; GJ\tilde f]
$,
thus $\det G=\det [g \; \;GJ\tilde f]$.
On the other hand,
$\det [g \;\; J\tilde g]=1$, so we also have
$\det G= \det [g \; \;\det G.J\tilde g]$. It follows that $\det [g \quad (GJ\tilde f-\det G. J\tilde g)] = 0$,
and therefore, for some $\beta \in \F$,
\[
GJ\tilde f = \beta g + \det G. J\tilde g.
\]
Multiplying this equation on the left by $\tilde g^T$, we get
$\beta= \tilde g^T GJ \tilde f$, since $\tilde g^T J \tilde g=0$.
\endpf

Note that any function $g$ belonging to $(H_\infty^\pm)^2$ or to $(H_2^\pm)^2$
is left invertible in $\F^{2 \times 1}$ if it is not identically zero.
Indeed if, for instance, the first component $g_{1\pm}$ of $g_\pm$ is not
identically zero, then we can take $\tilde g_\pm=(g^{-1}_{1+},0)$.

\begin{ex}\label{ex:3.3r}
{\rm
Let $\ds G=\begin{bmatrix}
\vphantom{\frac{|}{|}} 
\overline\theta & 0 \\ h & \overline r \end{bmatrix}$,
where $h \in L_\infty$, $\theta$ is an inner function, and $r(\xi)=\dfrac{\xi-i}{\xi+i}$ for $\xi\in \RR$.
Note that in this case the first component of any element in $\ker T_G$ belongs to the model space $K_\theta$.

We have $Gf=g$ with $f=(\theta,-h\theta r)$ and $g=(1,0)$
 and we can take
as their left inverses   the functions $\tilde f^T$ and $\tilde g^T$ with $\tilde f=(\overline\theta,0)$ and $\tilde g=(1,0)$.
We shall now use Theorem~\ref{thma:3.1} to describe $\ker T_G$. We have
\[
\SS= \{P^-(h\psi_+)+ \frac{k}{\xi-i}: \psi_+ \in H_2^+, k \in \CC \} \subset H_2^- 
\]
because from \eqref{eqa:3.2} we have, for $\psi_{1\pm}, \psi_{2\pm} \in H_2^\pm$,
\begin{eqnarray}  
&&\overline r\overline\theta [\theta \  -h\theta r] 
\left[ \begin{matrix} \psi_{1+} \\ \psi_{2+} \end{matrix} \right] = [1 \ 0]
\left[ \begin{matrix} \psi_{1-} \\ \psi_{2-} \end{matrix} \right]  \nonumber
\\
&\iff& \overline r \psi_{1+} - h \psi_{2+} = \psi_{1-} \nonumber \\
&\iff& \overline r \psi_{1+} - 2i \frac{\psi_{1+}(i)}{\xi-i} - P^+(h\psi_{2+})
= \psi_{1-}-2i \frac{\psi_{1+}(i)}{\xi-i} + P^-(h \psi_{2+}). \nonumber  \\ \label{eqny:*}
\end{eqnarray}
Since the left-hand side of this equation is in $H_2^+$ while the left-hand side is in $H_2^-$, both sides must be
equal to $0$, so we have from the right-hand side of \eqref{eqny:*} that
\beq\label{eqny:**}
\psi_{1-}= 
P^-(h\psi_{2+})+ \frac{k}{\xi-i} \with \psi_{2+} \in H_2^+, k \in \CC.
\eeq
Conversely, if $\psi_{1-}$ takes the form \eqref{eqny:**}, then $\psi_{1-} \in \SS$ because 
\eqref{eqny:*} holds with $\psi_{1+}= \dfrac{k}{\xi+i} + r P^+(h\psi_{2+})$.

We see that $\SS \ne \{0\}$, since $\dfrac{1}{\xi-i} \in \SS$. Given any $s \in \SS$ of the form given by the right-hand side of 
\eqref{eqny:**}, we have
\begin{eqnarray*}
\HH_+^s &=& \{(\phi_{1+},\phi_{2+}) \in (H_2^+)^2: -hr\phi_{1+}+sr  \in H_2^+\} \\
&=& \{(\phi_{1+},\phi_{2+}) \in (H_2^+)^2: P^-(hr\phi_{1+})=P^-(sr)\} 
\end{eqnarray*}
and $\HH_-^s = (H_2^-)^2$.
Since in this case $\tilde g^T(GJ\tilde f)=0$ we have from \eqref {eqa:3.6}   that
$\Lambda \in \tilde f^T \HH_+^s \cap \tilde g^T \HH_-^s$, which is equivalent to
\[
\Lambda=\overline\theta \phi_{1+} \with \phi_{1+} \in K_\theta, \quad  P^-(hr\phi_{1+})=P^-(sr).
\]
Therefore, from Theorem~\ref{thma:3.1}, specifically \eqref{eqb:3.5},
\[
\phi_+= \left(\phi_{1+}, \ -P^+(hr\phi_{1+})+ \frac{\tilde k}{\xi+i}\right),
\]
with $\phi_{1+} \in K_\theta$ and $\tilde k \in \mathbb C$, where we took into account the fact that 
$P^-(sr)=P^-(hr\phi_{1+})$ and $P^+(sr)=P^+(rP^-(h \psi_+))+\frac{k}{\xi+i}=\frac{\tilde k}{\xi+i}$ with $\tilde k \in \mathbb C$. Thus we have
\[
\ker T_G= \{  (\phi_{1+}, \ -P^+(hr\phi_{1+})): \phi_{1+} \in K_\theta \} +
\spam\left\{\left(0, \dfrac{1}{\xi+i}\right)\right\}.
\]
If $\theta$ is a finite Blaschke product of degree $n$, then $\dim\ker T_G=n+1$; otherwise $\dim\ker T_G=\infty$.

}\end{ex}

As a consequence of Theorem \ref{thma:3.1}, we have the following.

\begin{cor}\label{cora:3.3}
If
$
Gf=g$ and
\beq\label{eqa:3.16}
\left(\det G. f^T (H_2^+)^2 \right) \cap \left(g^T (H_2^-)^2 \right)= \{0\},
\eeq
and we define
\begin{eqnarray}
\HH_+ &=& \{ \psi_+ \in (H_2^+)^2: f \tilde f^T \psi_+ \in (H_2^+)^2 \}, \label{eqa:3.17}\\
\HH_- &=& \{ \psi_- \in (H_2^-)^2: g \tilde g^T \psi_- \in (H_2^-)^2 \}, \label{eqa:3.18} \\
\K &=& \tilde f^T \HH_+ \cap \tilde g^T \HH_-,\label{eqa:3.19}
\end{eqnarray}
then
\beq\label{eqa:3.20}
\ker T_G = \K f.
\eeq
\end{cor}
Since $0 \in \SS$, with $\SS$ defined in \eqref{eqa:3.2}, we also have the following consequence of Theorem \ref{thma:3.1}, which can be understood as
establishing a lower bound for $\ker T_G$.

\begin{cor}\label{cora:3.4}
If $Gf=g$ then, with the same notation as in Corollary \ref{cora:3.3}, we have
\[
\K f \subset \ker T_G.
\]
\end{cor}

By Coburn's Lemma \cite{coburn}, for any \T operator with scalar symbol $\phi \in L_\infty$, either $\ker T_\phi$ or
$\ker T_{\phi}^*=\ker T_{\overline\phi}$ is zero. It is well known that this property no longer holds when we consider \T operators with matrix symbol, since $T_G$ and $T_G^*=T_{\bar G^T}$ may both have a non-zero kernel.
However, using the result of Corollary \ref{cora:3.3}, we can state what may be seen as a version of Coburn's Lemma
for $2 \times 2$ block \T operators with symbol $G$. 
We shall need the following, which can easily be verified:

\begin{lem}\label{lema:3.5}
Let $G$ be a $2 \times 2$ matrix. Then
\beq
\det G. I = - GJG^T J.\
\eeq
\end{lem}

\begin{thm}\label{theor3.7}
Let $\det G \in \F\setminus \{0\}$. Then either $\ker T_G$ or $\ker T_{\overline G^T} $ 
is of scalar type.
\end{thm}
\beginpf
Assume that $\ker T_{\overline G^T} \ne \{0\}$ and let $\psi_+ \in \ker T_{\overline G^T}$, $\psi_+ \ne 0$. Then 
 we have
$\overline{G}^T\psi_+=\psi_- \in (H_2^-)^2$ and
\begin{eqnarray}\label{eqa:3.22}
\overline G^T \psi_+=\psi_-& \iff& G^T \,\overline{\psi_+}=\overline{\psi_-}
\iff GJG^TJ(J \overline {\psi_+})= -GJ \overline{\psi_-} \nonumber \\
&\iff& \det G. (J\overline{\psi_+})=GJ\overline{\psi_-}.
\end{eqnarray}
Therefore, $ GF_+= \det G. F_-$
with $F_\pm= J \overline{\psi_\mp} \in (H_2^\pm)^2$.
For any $\phi_+ \in \ker T_G$, we have
$G\phi_+=\phi_- \in (H_2^-)^2$, so
$
G[\phi_+ \;\; F_+]=[\phi_- \;\;\det G. F_-] $,
and it follows that 
\[
\det G. \det [\phi_+ \;\; F_+] = \det G. \det [\phi_- \;\;   F_-],
\]
i.e., on a set of positive measure in $\mathbb R$,
\beq\label{eqa:3.23}
\det [\phi_+ \;\; F_+] =   \det [\phi_- \;\;   F_-].
\eeq
Since the left-hand side of \eqref{eqa:3.23}
represents a function in $H_1^+$ while the right-hand side
represents a function in $H_1^-$, both are equal to zero.
Since $F_+$ and $F_-=GF_+$ admit left inverses because neither is identically equal to zero, it follows that every $\phi_+ \in \ker T_G$ is a scalar multiple of $F_+$.
\endpf

\begin{cor}\label{cora:3.7}
For every $G \in (L_\infty^{2 \times 2})$ with $\det G \in L_\infty \setminus \{0\}$, either $\ker T_G=\{0\}$,
or $\ker T_G^*=\{0\}$, or both kernels are of scalar type.
\end{cor}

\begin{cor}\label{cora:3.8}
If $\det G$ admits a (canonical) bounded factorization \cite{MP} $\det G=d_- d_+$ with $d_\pm \in \G H_\infty^\pm$, then both $\ker T_G$ and $\ker T_{\overline G^T}$ are of scalar type. In particular, $\ker T_G$ and $\ker T_{\overline G^T}$ are of scalar type whenever $\det G=1$.
\end{cor}
\beginpf
From \eqref{eqa:3.22}  we have that $\overline G^T \psi_+ = \psi_- \iff G(d_+^{-1} J \overline{\psi_-}) = d_- J \overline{\psi_+}$.
Since any $\phi_+ \in (H_2^+)^2$ can be written in the form
$\phi_+= d_+^{-1}J \overline{\psi_-}$ for some $\psi_- \in (H_2^-)^2$, and any $\phi_- \in (H_2^-)^2$  can be written
in the form $\phi_- = d_- J \overline{\psi_+}$ for some $\psi_+ \in (H_2^+)^2$, it follows that
$\ker T_G = \{0\}$ if and only if $\ker T_{\overline G^T}=\{0\}$. The result now follows from
Corollary \ref{cora:3.7}.
\endpf


Note that, for block Toeplitz operators, it is not always the case that a non-zero kernel can be given by a  symbol of determinant $1$,
as  the following simple example shows: let
$\ds G=\begin{bmatrix}
\overline r & 0 \\
\overline r & 0
\end{bmatrix}$.
The kernel of $T_G$ is $(k, h)$ where $k \in K_r$ and $h \in H_2^+$. The symbol can only have rows of the form
$(p \quad 0)$, for if $p,q \in L_\infty$ and $pk+qh \in H_2^-$ for all $k \in K_r$ and $h \in H_2^+$, then, taking $k=0$ we
see that 
$q=0$.
\\

While Corollary \ref{cora:3.3} provides sufficient conditions for the kernel of a Toeplitz operator with $2 \times 2$ symbol to be of scalar type, condition \eqref{eqa:3.16} is not a necessary one. To see this, let us consider the solution to $Gf=g$ that we obtain from \eqref{eq:2.1}
if we take the second columns of the matrix functions on the left and on the right hand sides of \eqref{eq:2.3}, instead of the first columns as was done in Section \ref{sec:2}. We get, using the same notation,
\[
Gf=g, \qquad \hbox{with} \qquad f= \begin{bmatrix} g_{12}^+ \\ g_{22}^+ \end{bmatrix}, \qquad 
g= r^k \begin{bmatrix} g_{12}^- \\ g_{22}^- \end{bmatrix}.
\]
Assuming, for simplicity, that $\det G=1$, we can choose $G_\pm$ such that $\det G_\pm=1$, and thus 
as left inverses for $f$ and $g$
we can take 
$\tilde f^T$ and $\tilde g^T$ given by
\[
\tilde f = J \begin{bmatrix} g_{11}^+ \\ g_{21}^+ \end{bmatrix}
= \begin{bmatrix} -g_{21}^+ \\ g_{11}^+ \end{bmatrix}, \qquad
\tilde g= r^{-k}J \begin{bmatrix} g_{11}^- \\ g_{21}^- \end{bmatrix}
= r^{-k} \begin{bmatrix} -g_{21}^- \\ g_{11}^- \end{bmatrix}.
\]

Applying Theorem \ref{thma:3.1}, we have $\SS = f^T (H_2^+)^2 \cap g^T (H_2^-)^2 = K_{r^k} \ne \{0\}$;
for each $s \in \SS=K_{r^k}$, we have 
$
\HH_+^s = (H_2^+)^2$ and
\[
 \HH_-^s = \left\{ \psi_- \in (H_2^-)^2: sr^{-k}  \begin{bmatrix} -g_{21}^- \\ g_{11}^- \end{bmatrix} \in (H_2^-)^2 \right \}
 = (H_2^-)^2,
 \]
 because $sr^{-k} \in H_\infty^-$ for $s \in K_{r^k}$. Therefore, from \eqref{eqa:3.6},
 \beq\label{eqa:3.24}
 \Lambda \in \tilde f^T (H_2^+)^2 = H_2^+,
 \eeq
 and, since in this case $GJ \tilde f = J \tilde g$, which implies that $\tilde g^TGJ \tilde f = 0$, we must also have
 \beq\label{eqa:3.25}
 \Lambda \in \tilde g^T (H_2^-)^2 = H_2^-.
 \eeq
 From \eqref{eqa:3.24} and \eqref{eqa:3.25} we get $\Lambda=0$ and it follows that $\ker T_G = K_{r^k} J \tilde f = K_{r^k} ( g_{11}^+\,,\, g_{21}^+)$
 as in Section \ref{sec:2}.
 
 The next result shows that every scalar-type \T kernel, for a $2 \times 2$ matrix symbol $G$, is of the form \eqref{eqa:3.20} with $\K$ given by \eqref {eqa:3.17}--\eqref {eqa:3.19},
 if $f$ and $g=Gf$ have left inverses.
 
 \begin{thm}\label{thma:3.9}
 If $\ker T_G = \K f$, where $f$ is a fixed function in $\F^{2 \times 1}$ such that $f$ and $g=Gf$ possess left
 inverses $\tilde f^T$ and $\tilde g^T$, respectively, and $\K \subset \F$, then
 \[
 \K = \tilde f^T \HH_+ \cap \tilde g^T \HH_-,
 \]
 where $\HH_\pm$ are defined as in \eqref{eqa:3.17}--\eqref{eqa:3.18}.
 \end{thm}
\beginpf
Let $k$ be any element of $\K$. Then $kf \in \ker T_G$ and we have
\beq\label{eqa:3.26}
kf = \psi_+ \in (H_2^+)^2, \qquad G(kf)= \psi_- \in (H_2^-)^2.
\eeq
From the first equation we get that $k=\tilde f^T \psi_+$, so
$f\tilde f^T \psi_+ = fk = \psi_+ \in (H_2^+)^2$; therefore $\psi_+ \in \HH_+$.

Analogously, from the second equation in \eqref{eqa:3.26}, we have $kg=\psi_- \in (H_2^-)^2$.
Therefore $k=\tilde g^T \psi_-$ and $\psi_-$ is such that $g\tilde g^T \psi_- = gk = \psi_- \in (H_2^-)^2$; thus $\psi_- \in \HH_-$.
We conclude that $\K \subset \tilde f^T \HH_+ \cap \tilde g^T \HH_-$.

Conversely, if $k \in f^T \HH_+ \cap \tilde g^T \HH_-$, then $kf \in \ker T_G$ (as in the last part of the proof
of Theorem \ref{thma:3.1}, with $s=0$) and, since $\ker T_G = \K f$, we have $kf= k_0 f$ with $ k_0 \in \K$. Multiplying 
on the left by
$\tilde f^T$ , we conclude that $k= k_0 \in \K$, so 
$f^T \HH_+ \cap \tilde g^T \HH_- \subset \K$.
\endpf

Naturally, one can say more about the space $\K$ if further assumptions are made on $f$ and $g$ in \eqref{eqa:3.1}.

\begin{thm}\label{thma:3.10}
If $f=\theta f_+$, where $\theta$ is an inner function and $f_+ \in (H_\infty^+)^2$, and $g=f_- \in (H_\infty^-)^2$, where
$f_\pm$ possess left inverses $\tilde f^T_\pm$ with $\tilde f_\pm \in (H_\infty^\pm)^2$, then
\[
\ker T_G = K_\theta f_+.
\]
\end{thm}
\beginpf
In this case we have $\HH_\pm=(H_2^\pm)^2$ and, if $\tilde f_\pm^T$ are
left inverses for $f_\pm$, then $\tilde f = \overline\theta \tilde f_+$ and $\tilde g=\tilde f_-$ provide
left inverses for $f$ and $g$ respectively; the result now follows from
\eqref{eqa:3.19} and \eqref{eqa:3.20}.
\endpf

We say that $f_\pm=(f_{1\pm},f_{2\pm}) \in (H_\infty^\pm)^2$ is a {\em corona pair\/} in $\CC^\pm$ if and only if
there exists $\tilde f_\pm \in (H_\infty^\pm)^2$ such that
$\tilde f^T_\pm f_\pm = 1$. In this case we say that
$f_\pm \in \CP^\pm$.
By the Corona Theorem, $f_\pm \in \CP^\pm$ if and only if
\beq\label{eq:3.27}
\inf_{z \in \CC^\pm} \left( |f_1^\pm(z)| + |f_2^\pm(z)| \right) > 0.
\eeq
Thus, under the conditions of Theorem \ref{thma:3.10}, we have $f_\pm \in \CP^\pm$.

The next theorem generalises Theorem \ref{thma:3.10}, establishing sufficient conditions 
for $\K$, in Corollary \ref{cora:3.3}, to be a model space or a shifted model space.
We shall use the following well-known result,
which 
 follows easily from the observation that if $\alpha,\beta$ are coprime inner functions
then $\alpha \phi_+ \in \beta H_2^+$ if and only if $ \phi_+ \in \beta H_2^+$ 
(a consequence of the uniqueness of the inner--outer factorization).

\begin{lem}\label{lema:3.11}
If $\phi_+ \in H_2^+$ and $\alpha_1$ and $\alpha_2$ are inner functions, then $\alpha_1 \phi_+ \in \alpha_2 H_2^+$  if and only 
if
$\phi_+ \in \ds \frac{\alpha_2}{\gamma_\alpha} H_2^+$, where
$\gamma_\alpha = \gcd\{ \alpha_1, \alpha_2 \}$.
\end{lem}

\begin{thm}\label{thma:3.12}
Let $\det G$ have a canonical bounded factorization $\det G=d_- d_+^{-1}$ (as in Section~\ref{sec:2}), and let $Gf=g$ with
componentwise
inner--outer factorizations
\beq
f=(\alpha_1 f_{1+}, \alpha_2 f_{2+}) , \qquad g=(\overline{\beta_1}f_{1-}, \overline{\beta_2}f_{2-}),
\eeq
where $\alpha_1, \alpha_2, \beta_1$ and $\beta_2$ are inner functions, $f_{1\pm},f_{2\pm}$ are outer functions in $H_\infty^\pm$, and
$(f_{1\pm},f_{2\pm}) \in \CP^\pm$. Then
\beq
\ker T_G = K_{\gamma_\alpha \gamma_\beta} 
\begin{bmatrix}
\ds \frac{\alpha_1}{\gamma_\alpha} f_{1+} \\
\ds \frac{\alpha_2}{\gamma_\alpha} f_{2+}
\end{bmatrix} =  K_{\gamma_\alpha \gamma_\beta} \,\overline{\gamma_\alpha}  f,
\eeq
where
\beq
\left \{ \begin{array}{rcl}
\gamma_\alpha &=& \gcd \{ \alpha_1, \alpha_2\} ,\\
\gamma_\beta &=& \gcd \{ \beta_1, \beta_2 \}  .
\end{array}
\right. 
\label{eqa:3.30}
\eeq
\end{thm}

\beginpf
In this case we have $\SS=\{0\}$ and $\ker T_G$ is given by Corollary \ref{cora:3.3}.
Let $[\tilde f_{1\pm} \; \tilde f_{2\pm}]$, with $(\tilde f_{1\pm} , \tilde f_{2\pm}) \in (H_\infty^\pm)^2$, be left inverses for
$(f_{1\pm}, f_{2 \pm})$, respectively, so that
$\tilde f=(\overline{\alpha_1}\tilde f_{1+}, \overline{\alpha_2}\tilde f_{2+})$ and
$\tilde g=( {\beta_1}\tilde f_{1-},  {\beta_2}\tilde f_{2-})$ are left inverses for $f$ and $g$, respectively.

For any $\psi_+ = (\psi_{1+},\psi_{2+}) \in (H_2^+)^2$ we have then
\begin{eqnarray*}
f \tilde f^T \psi_+ &=&
\begin{bmatrix}
f_{1+} \tilde f_{1+} \psi_{1+} + \alpha_1 \overline{\alpha_2} \tilde f_{2+} f_{1+} \psi_{2+} \\
\overline{\alpha_1}\alpha_2 \tilde f_{1+} f_{2+} \psi_{1+} + f_{2+} \tilde f_{2+} \psi_{2+}
\end{bmatrix}   \\
\end{eqnarray*}
so $f \tilde f^T \psi_+ \in (H_2^+)^2$ if and only if
\begin{eqnarray} 
&&\left\{ 
\begin{array}{rcl}
\overline{\alpha_2}\alpha_1 \tilde f_{2+}  f_{1+} \psi_{2+} &\in& H_2^+ ,\\
\overline{\alpha_1}\alpha_2 \tilde f_{1+}  f_{2+} \psi_{1+} & \in & H_2^+,
\end{array}
\right. \nonumber \\
&\iff&
\left\{ 
\begin{array}{rcl}
\alpha_1 \tilde f_{2+}  f_{1+} \psi_{2+} &\in& \alpha_2 H_2^+, \\
 \alpha_2\tilde f_{1+} f_{2+} \psi_{1+}  & \in & \alpha_1 H_2^+,
\end{array}
\right. \nonumber\\
&\iff&
\left\{ 
\begin{array}{rcl}
 \tilde f_{2+}  f_{1+} \psi_{2+} &\in& \frac{\alpha_2}{\gamma_\alpha} H_2^+, \\
 \tilde f_{1+}   f_{2+} \psi_{1+}  & \in & \frac{\alpha_1}{\gamma_\alpha} H_2^+,
\end{array}
\right.\label{eqa:3.31}
\end{eqnarray}
with $\gamma_\alpha$ defined by \eqref{eqa:3.30}. Thus $(\psi_{1+},\psi_{2+}) \in \HH_+$ if and only if $\psi_{1+}, \psi_{2+} \in H_2^+$ and 
\[
\psi_{1+} = \frac{\alpha_1}{\gamma_\alpha}\frac{\phi_{1+}}{\tilde f_{1+} f_{2+}}, \qquad  \psi_{2+} = \frac{\alpha_2}{\gamma_\alpha} \frac{\phi_{2+}}{\tilde f_{2+}f_{1+}},
\qquad {\rm with} \quad \phi_{1+}, \phi_{2+} \in H_2^+, 
\]
where, since
$\frac{\phi_{1+}}{ f_{2+}}=\psi_{1+} \frac{\gamma_\alpha}{\alpha_1}\tilde f_{1+}\,\in\mathcal N^+\cap L_2$
and 
$\frac{\phi_{2+}}{ f_{1+}}=\psi_{2+} \frac{\gamma_\alpha}{\alpha_2}\tilde f_{2+}\,\in\mathcal N^+\cap L_2,$
we have
\beq\label{eqa:A}
\frac{\phi_{1+}}{ f_{2+}}\,,\,\frac{\phi_{2+}}{ f_{1+}}\,\in H_2^+.
\eeq
Analogously, we get that  $(\psi_{1-},\psi_{2-}) \in \HH_-$ if and only if $\psi_{1-}, \psi_{2-} \in H_2^-$ with
\[
\psi_{1-} = \overline{\left( \frac{\beta_1}{\gamma_\beta}\right)}\frac{\phi_{1-}}{\tilde f_{1-}f_{2-}}\,, \qquad  \psi_{2+} = \overline{\left(\frac{\beta_2}{\gamma_\beta} \right)}\frac{\phi_{2-}}{\tilde f_{2-}f_{1-}}\,,
\qquad  \phi_{1-}, \phi_{2-} \in H_2^-,
\]
where
\beq\label{eqa:AA}
\frac{\phi_{1-}}{ f_{2-}}\,,\,\frac{\phi_{2-}}{ f_{1-}}\,\in H_2^-.
\eeq
Thus $\K= \tilde f^T \HH_+ \cap \tilde g^T \HH_-$ consists of the functions $k$ such that
\beq\label{eqa:3.33}
k= [\overline{\alpha_1}\tilde f_{1+} \quad \overline{\alpha_2} \tilde f_{2+}]
\begin{bmatrix}
\frac{\alpha_1}{\gamma_\alpha}\frac{\phi_{1+}}{\tilde f_{1+} f_{2+}} \\\frac{\alpha_2}{\gamma_\alpha}\frac{\phi_{2+}}{\tilde f_{2+} f_{1+}}
\end{bmatrix}
=
[ \beta_1 \tilde f_{1-} \quad \beta_2 \tilde f_{2-} ]
\begin{bmatrix}
 \overline{\left( \frac{\beta_1}{\gamma_\beta}\right)}\frac{\phi_{1-}}{\tilde f_{1-}f_{2-}} \\  \overline{\left( \frac{\beta_2}{\gamma_\beta}\right)}\frac{\phi_{2-}}{\tilde f_{2-}f_{1-}}
\end{bmatrix},
\eeq
i.e.,
\[
k= \overline{\gamma_\alpha} \left(\frac{\phi_{1+}}{f_{2+}}+ \frac{\phi_{2+}}{f_{1+}}\right) = \gamma_\beta  \left(\frac{\phi_{1-}}{f_{2-}}+\frac{\phi_{2-}}{f_{1-}}\right).
\]
Taking \eqref{eqa:A} and \eqref{eqa:AA} into account, it follows that
$\K \subset \overline{\gamma_\alpha}H_2^+ \cap \gamma_\beta H_2^-$.

Conversely, if $k=\overline{\gamma_\alpha}\phi_+ = \gamma_\beta \phi_-$, with $\phi_\pm \in H_2^\pm$, we can write
\[
k= \overline{\gamma_\alpha} \phi_+ =  \overline{\gamma_\alpha}\tilde f^T f  \phi_+ = \tilde f^T (\overline{\gamma_\alpha} f \phi_+),
\]
where 
$ \overline{\gamma_\alpha} f \phi_+ \in \HH_+$ because $\overline{\gamma_\alpha} f \phi_+ \in (H_2^+)^2$ and
$f \tilde f^T (\overline{\gamma_\alpha} f \phi_+) = f\overline{\gamma_\alpha}  \phi_+$.

Thus $k \in \tilde f^T \HH_+$ and,
analogously, we can show that $k \in \tilde g^T \HH_-$, so that
$k \in \tilde f^T \HH_+ \cap \tilde g \HH_- = \K$. We conclude that
$\K = \overline{\gamma_\alpha} H_2^+ \cap \gamma_\beta H_2^- = \overline{\gamma_\alpha} K_{\gamma_\alpha \gamma_\beta}$.
\endpf

\begin{cor}
With the same assumptions as in Theorem \ref{thma:3.12}, $T_G$ is injective if and only if $\gamma_\alpha$ and $\gamma_\beta$
are constant.
\end{cor}

\begin{cor}\label{cora:3.14}
If the assumptions of Theorem \ref{thma:3.12} are satisfied, with
\[
f=\begin{bmatrix}
\alpha \\ O_+  
\end{bmatrix}, \qquad 
g=\begin{bmatrix} 
\overline\beta  \\ O_-  
\end{bmatrix},
\]
where $\alpha,\beta$ are inner functions and $O_+, \overline{O_-} \in H_\infty^+$ are outer functions, then
$T_G$ is injective.
\end{cor}

Note that, if $\det G$ has a canonical bounded factorization and $f\in \CP^+$ and $g\in \CP^-$, then $T_G$ is invertible (\cite{CDR}). Roughly speaking, Corollary \ref{cora:3.14} means that if we have $Gf_+=f_-$ with $f_\pm \in (H_\infty^\pm)^2$, and  the two components of $f_\pm$  ``approach zero simultaneously at some point'' in $\CC^\pm \cup \RR$ (so that they do not satisfy the corona conditions \eqref{eq:3.27}),we may still have an injective \T operator
as long as they do not ``approach zero simultaneously'' through a common inner factor.\\

Suppose now that $\ker T_G$ is of scalar type, $\ker T_G=\K f $ as in Theorem \ref{thma:3.9}.
If $\ker T_G \ne \{0\}$, we may now ask, in view of the results of Theorems \ref{thma:3.10} and \ref{thma:3.12}, whether a scalar type $\ker T_G$ can always
be described
as the product of some fixed $2 \times 1$ function by a scalar \T kernel, and in particular a model space.

Toeplitz kernels constitute an important subset of the class of nearly $S^*$-invariant subspaces of $H^+_2$. Here by $S^*$ we denote the backward shift operator $S^*=P^+ \bar r P^+_{|H^+_2}$; a  subspace $M$ of $H^+_2$ is nearly $S^*$-invariant if and only if $S^* \phi_+ \in M$ for all $\phi_+ \in M$ such that $\phi_+(i)=0$. 
Hitt proved (in the unit disk setting) that any nontrivial closed nearly $S^*$-invariant subspace of $H^+_2$ has the form $M=hK_\theta$ where $\theta$ is an inner function vanishing at $i$, $\frac{h}{\xi - i}\in M$ is the element of unit norm with positive value at $i$ which is orthogonal to all elements of $M$ vanishing at $i$, and $h$ is an isometric multiplier from $K_\theta$ into $H^+_2$ (\cite {Hitt,Sarason}).

In the next theorem we show that, if $\ker T_G$ is of scalar type, then it is the product of a scalar nearly $S^*$-invariant subspace (which can be explicitly described) with a fixed $2 \times 1$ function.

Note that, when $\ker T_G$ is of scalar type and $\ker T_G \ne \{0\}$, then there exist $F_\pm \in (H_2^\pm)^2 \setminus \{0\}$
such that $GF_+=F_-$ and, if $F_\pm=(F_{1\pm},F_{2\pm})$, then either
$F_{1\pm}$ or $F_{2\pm}$ is not identically zero, and thus $F_\pm$ has a left inverse defined a.e.\ on $\RR$.
We shall assume, without loss of generality, that $F_{1\pm} \ne 0$, so that $F^{-1}_{1\pm}$ are defined a.e.\ on $\RR$.
If $\ker T_G=\K f$ with $f=(f_1,f_2)$ and $Gf=g=(g_1,g_2)$, as in Theorem \ref{thma:3.9}, then
\[
F_{j+}=k_0 f_j, \qquad F_{j-}=k_0 g_j \qquad (j=1,2),
\]
where $k_0 \in \K$ and $k_0,f_1,g_1$ are different from zero a.e.\ on $\RR$. 

We shall also use the following notation:
if $F_\pm \in H_2^\pm \setminus \{0\}$, then we write $F_\pm = I_\pm O_\pm$, where
$I_+, \overline{I_-}$ are inner functions and $O_+, \overline{O_-}$ are outer functions in $H_2^+$.
If $F_\pm=0$ then we write $F_\pm=I_\pm O_\pm$ with $I_\pm,O_\pm=0$.
If $\alpha$ is an inner function, then $\gcd\{\alpha,0\}=\alpha$.


\begin{thm}\label{thma:3.15}
Let 
 $\ker T_G$ be of scalar type, 
$\ker T_G \ne \{0\}$, with
$\ker T_G=\K f$ as in Theorem \ref{thma:3.9}.
Then there exist an $F \in \F^{2 \times 1}$ and a nearly $S^*$-invariant subspace $\tilde \K \subset H_2^+$, which is closed if $G\in L_\infty^{2\times 2}$, such that
$\ker T_G= \tilde \K F$.
Moreover, if $F_+=(F_{1+},F_{2+}) \ne 0$ is
a given element of $\ker T_G$ and $GF_+=F_-=(F_{1-},F_{2-})$ with
$F_{j\pm}=I_{j\pm}O_{j\pm}$ ($j=1,2$), using the notation above, and we suppose that $F_{1\pm} \ne 0$, then
\beq\label{eqa:3.35a}
\ker T_G= \tilde \K \frac{F_+}{\gamma_+ O_{1+}}= \tilde \K \, \left(\frac{I_{1+}}{\gamma_+},\frac{F_{2+}}{\gamma_+ O_{1+}}\right),
\eeq
where
\beq\label{eqa:3.35}
\tilde \K= \left\{ \tilde\psi_+ \in \ker T_{\gamma_- \overline{\gamma_+} O_{1-}/O_{1+}}
\cap
\ker T_{\gamma_- \overline{\gamma_+} O_{2-}/O_{1+}}: \frac{O_{2+}}{O_{1+}} \tilde \psi_+ \in H_2^+ \right \},
\eeq
\beq\label{eqa:3.36}
\gamma_+ = \gcd(I_{1+},I_{2+}), \qquad \overline{\gamma_-}=\gcd(\overline{I_{1-}},\overline{I_{2-}}).
\eeq
\end{thm}

\beginpf
Let $F_+=(F_{1+},F_{2+}) \ne 0$ belong to $\ker T_G$ and let 
$GF_+=F_-=(F_{1-},F_{2-})$, where we assume that $F_{1\pm} \ne 0$.
Since we can write $f=k_0^{-1}F_+$ with $k_0 \in \K$, we have
$\ker T_G=(\K k_0^{-1})F_+$, and we can then assume that $\ker T_G=\K F_+$ and apply Theorem
\ref{thma:3.9}.

Defining $\HH_\pm$ as in \eqref{eqa:3.17}--\eqref{eqa:3.18}, we have, for $f=F_+$ and $\tilde f=(F_{1+}^{-1},0)$,
\begin{eqnarray*}
\psi_+ \in \HH_+ & \iff & f \tilde f^T \psi_+ \in (H_2^+)^2, \quad \psi_+ \in (H_2^+)^2 \\
&\iff& (\tilde f^T \psi_+)f \in (H_2^+)^2, \quad \psi_+ \in (H_2^+)^2 \\
&\iff& F_{1+}^{-1} \psi_{1+} (F_{1+},F_{2+}) \in (H_2^+)^2, \quad \psi_{1+} \in H_2^+, \quad \psi_{2+}\in H_2^+ \\
& \iff & \psi_{1+}F_{2+} \in F_{1+}H_2^+, \quad  \psi_{1+} \in H_2^+, \quad \psi_{2+}\in H_2^+  .
\end{eqnarray*}
We have
\beq\label{eq:eqa:3.37}
\psi_{1+}F_{2+} \in F_{1+}H_2^+ \iff \psi_{1+} I_{2+} \frac{O_{2+}}{O_{1+}} \in I_{1+}H_2^+.
\eeq
In this case $\psi_{1+}O_{2+}/O_{1+} \in \N_+ \cap L_2 = H_2^+$ and it follows from the second relation
in \eqref{eq:eqa:3.37}
and Lemma \ref{lema:3.11} that the right hand side of \eqref {eq:eqa:3.37} is equivalent to
\[
\psi_{1+} \frac{O_{2+}}{O_{1+}}= \frac{I_{1+}}{\gamma_+} \phi_+ \quad \hbox{with} \quad \phi_+ \in H_2^+ \quad \hbox{and} \quad
\gamma_+ = \gcd\{ I_{1+},I_{2+} \}.
\]
Since $\psi_{1+} \in H_2^+$ and 
\[
\frac{O_{1+}}{O_{2+}}\frac{I_{1+}}{\gamma_+} \phi_+ \in H_2^+ \iff  
\frac{O_{1+}}{O_{2+}}  \phi_+ \in H_2^+
\]
we conclude that 
\beq
\begin{array}{rll}
& \psi_{1+}F_{2+} \in F_{1+}H_2^+, \quad & \psi_{1+} \in H_2^+ \\ \\
\iff & \ds \psi_{1+} = \tilde\psi_+ \frac{I_{1+}}{\gamma_+}, \quad \hbox{with} \quad & \tilde\psi_+ \in H_2^+, \quad \ds \frac{O_{2+}}{O_{1+}} \tilde \psi_+ \in H_2^+. \label{eqa:3.38}
\end{array}
\eeq
So, $\psi \in \tilde f^T\HH_+$ if and only if $\psi=F_{1+}^{-1} \psi_{1+}$, where $\psi_{1+}$
satisfies \eqref{eqa:3.38}, i.e., 
\beq  \label{3.41}
\tilde f^T \HH_+ = \left \{
\psi \in \F: \psi= \frac{\tilde\psi_+}{\gamma_+ O_{1+}} \;\;\hbox{with} \;\; \tilde\psi_+ \in H_2^+, \;\; \frac{O_{2+}}{O_{1+}} \tilde \psi_+ \in H_2^+
\right \}.
\eeq
Analogously we get, for $F_{2-} \ne 0$ and $\gamma_-$ defined in \eqref{eqa:3.36},
\beq \label{3.42}
\tilde g^T \HH_- = \left \{
\psi \in \F: \psi= \frac{\tilde\psi_-}{\gamma_- O_{1-}} \;\; \hbox{with} \;\; \tilde\psi_- \in H_2^-, \;\; \frac{O_{2-}}{O_{1-}} \tilde \psi_- \in H_2^-
\right \}.
\eeq
Therefore, for $\psi$ to belong to $\K=\tilde f^T \HH_+ \cap \tilde g^T \HH_-$, the functions $\tilde \psi_\pm \in H_2^\pm$ in \eqref {3.41}--\eqref {3.42} must be such that
\[
\frac{\tilde \psi_+}{\gamma_+ O_{1+}} = \frac{\tilde \psi_-}{\gamma_- O_{1-}}, 
\]
i.e., $\ds \tilde \psi_+ \in \ker T_{\gamma_- \overline{\gamma_+} O_{1-}/O_{1+}}$, and the condition
$\tilde\psi_- O_{2-}/O_{1-} \in H_2^-$ in \eqref {3.42} can be expressed by
\[
\gamma_- \overline{\gamma_+} \frac{O_{2-}}{O_{1+}} \tilde \psi_+ \in H_2^-, \quad \hbox{i.e.,} \quad \tilde \psi_+ \in \ker T_{{\gamma_-}\overline{\gamma_+} O_{2-}/O_{1+}}.
\]

Finally, taking into account the last condition in \eqref{eqa:3.38}, we have
\[
\K = \tilde f^T \HH_+ \cap \tilde g^T \HH_- = \frac{1}{\gamma_+ O_{1+}} \tilde \K,
\]
where, for $F_{2+}, F_{1+} \ne 0$,
\[
\tilde \K = \left\{
\tilde\psi_+ \in \ker T_{\gamma_- \overline{\gamma_+} O_{1-}/O_{1+}}
\cap
\ker T_{\gamma_- \overline{\gamma_+} O_{2-}/O_{1+}}:
\frac{O_{2+}}{O_{1+}} \tilde \psi_+ \in H_2^+
\right\}.
\]
It is easy to see that $\tilde \K$ is nearly $S^*$-invariant, because
Toeplitz kernels are nearly $S^*$-invariant subspaces and  if $\frac{O_{2+}}{O_{1+}}\tilde\psi \in H_2^+$ and $r^{-1} \tilde\psi_+ \in H_2^+$, then $\frac{O_{2+}}{O_{1+}}\tilde\psi r^{-1} \in \N^+ \cap L_2 = H_2^+$.

If $F_{2\pm}=0$, then $\HH_\pm = (H_2^\pm)^2$ and we again find that $\K=\tilde\K \ds \frac{F_+}{\gamma_+ O_{1+}}$,
where
\begin{eqnarray*}
\tilde\K &=& \ker T_{\gamma_- \overline{\gamma_+}O_{1-}/O_{1+}}\cap
\ker T_{\gamma_- \overline{\gamma_+} O_{2-}/O_{1+}}   \\
&=&\ker T_{\gamma_-  O_{1-}/F_{1+}}\cap \ker T_{\gamma_-   \frac{O_{2-}}{F_{1+}}} \quad \hbox{if} \quad F_{2+}=0, \quad F_{2-} \ne 0, \\
\noalign{ \hbox{and}} \\
\tilde \K &=& \ker  T_{\gamma_- \overline{\gamma_+}O_{1-}/O_{1+}} = \ker T_{F_{1-}/F_{1+}}
\quad \hbox{if} \quad F_{2+}=0, \quad F_{2-} = 0.
\end{eqnarray*}
Finally, if $G\in L_\infty^{2\times 2}$ then $\ker T_G$ is closed, and it follows from \eqref{eqa:3.35a} that $\tilde \K \, \frac{I_{1+}}{\gamma_+}$ is closed, so $\tilde \K$ is closed.
\endpf


Since $\tilde\K$ is nearly-invariant, it can be expressed in terms of a Toeplitz kernel by Hitt's
theorem, say $\tilde K=K g$, where $K$ is a Toeplitz kernel and $g$ is a scalar function (an isometric multiplier);
however, there is no reason to suppose that $\tilde K$ is already a Toeplitz kernel.

Related to these results, a very natural question regarding scalar type Toeplitz kernels is whether they have a maximal function. It was proved in \cite{CP14} that for every $\phi_+  \in H^+_2$ there exists a so called minimal kernel $K_m(\phi_+)$ such that very other kernel $K$ with $\phi_+ \in K$ contains $K_m(\phi_+)$. We say that $\phi_+$ is a maximal function for $K$ if $K=K_m(\phi_+)$; every scalar Toeplitz kernel has a maximal function. For scalar type Toeplitz kernels we have the following, taking the result of Theorem \ref{thma:3.15} into account.

\begin{thm}\label{thma:3.17}
If $\ker T_G$, with $G\in L_\infty^{2\times 2}$, is of scalar type, then there exists a maximal function for $\ker T_G$.
\end{thm}

\beginpf
Let $\ker T_G = \tilde\K F= Kf$, as above,
where $f=gF$ as above,
$K$ is a scalar Toeplitz kernel with maximal
function $\phi_+ $;
 let $\phi_+= IO$, where $I$ is inner and $O$ is an outer function in $H^+_2$. Then $K=\ker T_{(\bar I \bar O /O)}$ (\cite {CP14})
and every element in $K$ has the form $[IO/\bar O]\psi_-$, where
$\psi_- \in H_2^-$.

Obviously $\phi_+ f$ belongs to $\ker T_G$. Suppose that $\phi_+f = IOf$ belongs to the kernel of some $T_H$ with $H\in L_\infty^{2\times 2}$; then
$H( IOf ) = \phi_-\in (H_2^-)^2$.  We want to prove that $\ker T_G$ is contained in $\ker T_H$.
 Take any element  $[IO/\bar O]\psi_- f \in Kf$; we have $H
[IO/\bar O]\psi_- f= (\psi_-/ \bar O) H (IOf) = 
(\psi_-/ \bar O)  \phi_-$. Since this is (componentwise)  in the Smirnov class $\N_-$ and in
$L_2$, it is in  $(H_2^-)^2$. Therefore   $[IO/\bar O]\psi_- f$ is in
$\ker T_H$.

\endpf

Note that, as remarked by R.~O'Loughlin, not all block Toeplitz kernels have a maximal function. This is an immediate consequence of Theorem 5.5 and Corollary 5.3 in \cite{CP14}.


\section{Applications to  truncated Toeplitz\\ operators}

Let $h \in \F$ and, for any inner function $\theta$, let
\beq
\D_\theta= \{ f_\theta \in K_\theta: h f_\theta \in L_2 \}.
\eeq
The operator $A_h^\theta: \D_\theta \to K_\theta$ defined by
\beq
A_h^\theta f_\theta = P_\theta (h f_\theta), \quad f_\theta \in \D_\theta,
\eeq
where $P_\theta: L_2 \to K_\theta$ denotes the orthogonal projection, 
is called the {\em truncated Toeplitz operator\/} (in $K_\theta$) with symbol $h$. 
If $h$ belongs to the Sobolev space $\lambda_+ L_2$, where $\lambda_+(\xi)=\xi+i$, then
$A_h^\theta$ is densely defined on $K_\theta$; if $h \in L_\infty$, then $A_h^\theta$ is a 
bounded operator on $K_\theta$.

It is clear that $\phi_{1+} \in \ker A^\theta_h$ if and only if
$\phi_{1+} \in H_2^+$ and the following two conditions hold:
\beq
\left\{
\begin{array}{rcl}
\overline\theta\phi_{1+} &=& \phi_{1-}, \\
h{\phi_1+}&=& \phi_{2-}-\theta \phi_{2+}, \quad \hbox{with} \quad \phi_{1-},\phi_{2-} \in H_2^-, \quad \phi_{2+} \in H_2^+.
\end{array}
\right.
\eeq
Therefore, $\ker A^\theta_h$ consists of the first components of the elements in the kernel of the
Toeplitz operator $T_G$ with
\beq\label{eqa:4.4}
G= \begin{bmatrix}
\overline\theta & 0 \\ h & \theta \end{bmatrix},
\eeq
defined on $\D= \{\Phi_+ \in (H_2^+)^2: G\Phi_+ \in (L_2)^2 \}$.
In particular, we have that $\ker A^\theta_h=\{0\}$ if and only if
$\ker T_G=\{0\}$.
Thus we can apply the results of Section \ref{sec:3}, for $G$ of the form \eqref{eqa:4.4}, to study the kernels of truncated Toeplitz operators. Note that we have $Gf=g$ with
$f=(f_1,f_2)$ and $g=(g_1,g_2)$ if and only if 
$g_1=\overline\theta f_1$ and $h=\ds\frac{g_2-\theta f_2}{f_1}$.

The following is a consequence of Theorem \ref{thma:3.15} and Corollary \ref{cora:3.8}.

\begin{thm}\label{t4.1}
The kernel of any truncated \T operator is the product of a nearly 
$S^*$-invariant subspace of $H_2^+$, given in \eqref{eqa:3.35} and \eqref{eqa:3.36},  by an inner function.
\end{thm}
\beginpf
Let $h$ be the symbol of the truncated \T operator $A^\theta_h$, and let $G$ be defined by \eqref{eqa:4.4}. By Corollary \ref{cora:3.8}, $\ker T_G$ is of scalar type and,
if $\ker T_G\ne \{0\}$,
then by Theorem \ref{thma:3.15} we have $\ker T_G= \ds \tilde\K \frac{F_+}{\gamma_+ O_{1+}}$, where
$F_+ = (F_{1+},F_{2+})$ is a given function in $(H_2^+)^2$ with $F_{1+}=I_{1+}O_{1+} \in K_\theta\,, F_{1+}\ne 0, $
 $\gamma_+ \preceq I_{1+}$, and $\tilde\K$ is given by \eqref{eqa:3.35} and \eqref{eqa:3.36}.
Therefore 
$\ds\ker T_G= \tilde\K  \left( \frac{I_{1+}}{\gamma_+},\frac{F_{2+}}{\gamma_+O_{1+}}\right)$, and 
$\ds\ker A^\theta_h = \frac{I_{1+}}{\gamma_+}\tilde\K$.
\endpf

\begin{rem}{\rm
Recently, Ryan O'Loughlin \cite{ryan} has arrived at a similar result by a different route.
Namely, it follows from \cite[Cor. 4.5]{CCP} that the kernel of a $2 \times 2$ block Toeplitz
operator $T_G$ can be written as $F_0 ((H_2^+)^r \ominus \Theta  (H_2^+)^{r'})$, where
$r$ and $r'$ are integers with $1 \le r' \le r \le 2$,
$\Theta \in (H_\infty^+)^{r \times r'}$ is inner, and $F_0 \in (H_2^+)^{2 \times r}$; in the
case 
$\ds G=\begin{bmatrix}\overline\theta & 0 \\ g & \theta
\end{bmatrix}$, it is possible to take $r'=r=1$, although no explicit formula for $\Theta$ is given.\\
Note that if $\K_1 \phi_+ = \K_2 \psi+$, 
where $\K_1$ and $K_2$ are model spaces and $\phi_+=(\phi_{1+},\phi_{2+})$,
$\psi_+=(\psi_{1+},\psi_{2+})$ are analytic in $\CC^+$,
then 
we have multipliers 
$\psi_{+1}/\phi_{1+}$ and $\psi_{+2}/\phi_{2+}$ 
from one model space
onto another, so that the work in \cite{crofoot}, \cite{FHR} and \cite{CP18} can be applied. Indeed,
if $w K_\alpha=K_\beta$, then $\beta=c \alpha w/\overline w$, where $c$ is a unimodular constant.\\
} 
\end{rem}


Regarding the formula \eqref{eqa:3.35} for $\tilde\K$ in Theorem \ref{thma:3.15}, note that,
for $G$ of the form \eqref{eqa:4.4},  if $\ker T_G \ne \{0\}$ then $F_{1\pm} \ne 0$ and we have $\overline\theta F_{1+}=F_{1-} $ if and only if $\frac{O_{1-}}{O_{1+}}= \overline\theta I_{1+}\overline{I_{1-}}$.
Therefore the symbol $\ds \gamma_- \overline{\gamma_+} \frac{O_{1-}}{O_{1+}}$ in  \eqref{eqa:3.35} is bounded.
Moreover, from \eqref{eqa:3.35} we have:
\[
\tilde \K= \left\{ \tilde\psi_+ \in \ker T_{\overline\theta (I_{1+}/\gamma_+)(\overline{I_{1-}/\gamma_-})}\cap \ker T_{
 \gamma_- \overline{\gamma_+} \frac{O_{2-}}{O_{1+}}} \,:\, \frac{O_{2+}}{O_{1+}} \tilde \psi_+ \in H_2^+ \right\},
\]
which takes the form
\beq\label{4.4a}
\tilde \K=  \ker T_{\overline\theta (\overline{I_{1-}/\gamma_-})}\cap \ker T_{
\gamma_-  \frac{O_{2-}}{F_{1+}}} \quad\hbox{if} \quad F_{2+}=0, F_{2-}\ne0,
\eeq
and
\[
\tilde\K= \ker T_{\overline\theta}=K_\theta \quad \hbox{if} \quad F_{2+}=F_{2-} = 0.
\]

It may happen that the inner function mentioned in Theorem \ref{t4.1} is necessarily a constant, as it happens if the symbol $h$ of the truncated Toeplitz operator is in $H_\infty^-$. In that case, for $\gamma = \gcd(\theta,(\bar h)_i)$, it is easy to see that $F_+=(F_{1+},F_{2+})=(\frac {\gamma-\gamma(i)}{\xi-i},0) \in \ker T_G$, with $G$ given by \eqref{eqa:4.4},
and $GF_+=(F_{1-},F_{2-})=\frac {1-\bar\gamma\gamma(i)}{\xi-i}(\bar\theta\gamma,h\gamma)$. With the notation of Theorem \ref{thma:3.15}, we have $\gamma_-=1\,,\,I_{1-}=\bar \theta\gamma\,,\,O_{2-}=\frac {1-\bar\gamma\gamma(i)}{\xi-i}$, and it follows from \eqref{4.4a} that $\ker A_h^\theta=\ker T_{\bar \theta}=K_\gamma$
(as may also be verified by direct calculation).\\

We now apply the previous results to studying the kernels of some classes of TTO. Our motivation for the 
examples that we
shall consider is the following. The kernels of TTO with
so-called $\theta$-separated symbols, of the form
\beq
h= \overline\alpha h_1 + \beta h_2, \qquad (h_1 \in H_\infty^-, \quad h_2 \in H_\infty^+),
\eeq
where
\beq
\alpha\beta \succeq \theta,
\eeq
were studied in \cite{CP16}. We  consider here two cases where $h_1 \in H_\infty^+$ and $h_2 \in H_\infty^-$. In the
first case we assume that $h_1$ and $\overline{h_2}$ are inner, with 
$h_1 \prec \alpha$ and $\overline{h_2} \prec \beta$.
In the second case we assume that
$h_1$ and $h_2$ are rational functions, $h_1=R_1^+ \in  \R^+$, and 
$h_2=R_2^- \in \R^-$ with $\R^\pm=\R \cap H_\infty^\pm$, where $\R$ is the set of rational functions, which generalizes the study of truncated Toeplitz operators with $\theta$-separated symbols to the case where $h_1$ and $h_2$ admit poles in the lower and upper half planes, respectively.\\

\subsection{The first case: $A_h^\theta$ with $ h= C_1 \overline \alpha + C_2 \beta, \qquad $}\label{sec:4.1}

Toeplitz operators with almost-periodic symbols   of the form  \eqref{eqa:4.4}
where $\theta(\xi)=e^{i\lambda \xi}$, with $\lambda \in \RR$, and
\[
h(\xi)=C_1 \exp(-ia\xi)+C_2 \exp(ib \xi), \with a,b \in \RR^+ \quad \hbox{and} \quad a,b < \lambda,
\]
 have been studied by several authors (see, e.g., \cite{BKS,D MSc}). In this section we 
 generalize this class by studying symbols of the form  \eqref{eqa:4.4} with 
\[
h= C_1 \overline \alpha + C_2 \beta, \qquad (C_1,C_2 \in \CC \setminus \{0\}),
\]
where $\alpha$, $\beta$ are non-constant inner functions satisfying   the following conditions:
\beq\label{eq:43}
\alpha, \beta \prec \theta,
\eeq
\beq\label{eq:44}
\hbox{for some} \quad n \ge 1, \qquad (\alpha\beta)^{n-1} \prec \theta, \qquad (\alpha\beta)^n \succeq \theta,
\eeq
\beq\label{eq:45}
\hbox{either} \quad \alpha^n \beta^{n-1} \preceq \theta \qquad \hbox{or} \quad \alpha^n \beta^{n-1} \succeq \theta,
\eeq
\beq\label{eq:46}
\hbox{either} \quad \alpha^{n-1} \beta^{n} \preceq \theta \qquad \hbox{or} \quad \alpha^{n-1} \beta^{n} \succeq \theta.
\eeq
If $\gamma, \delta$ are inner functions such that either $\gamma \preceq \delta$ or $\delta \preceq \gamma$, we say that $\min \{\gamma,\delta\}=\gamma$ if $\gamma \preceq \delta$ and  $\min \{\gamma,\delta\}=\delta$ if $\delta \preceq \gamma$.
With this notation, if \eqref{eq:45} holds then either $\alpha^n \preceq \theta \overline\beta^{n-1}$ or
$\theta\overline\beta^{n-1} \preceq \alpha^n$, so in the first case 
$\min \{\theta\overline\beta^{n-1},\alpha^n\}=\alpha^n$, and in the
second case  $\min \{\theta\overline\beta^{n-1},\alpha^n\}=\theta\overline\beta^{n-1}$.

Analogously, if \eqref{eq:46} holds, then either $\alpha^{n-1}\beta^n \overline\theta \preceq \alpha$
or $\alpha \preceq \alpha^{n-1}\beta^n \overline\theta$, and 
$\min\{ \alpha^n\beta^n \overline\theta, \alpha \}$\,,\, $\min \{\beta, \theta \overline\alpha^{n-1}\overline\beta^{n-1} \}$ also exist.

Note that, if $\epsilon$ is a singular inner function (for instance, an exponential  $\exp (i\lambda\xi)$ for a given positive real $\lambda$) and $\lambda,a,b$ are real positive numbers with $a,b<\lambda$, then $\theta=\epsilon^\lambda\,,\,\alpha=\epsilon^a\,,\,\beta=\epsilon^b$ always satisfy \eqref{eq:43}-\eqref{eq:46} if we take $n$ to be the smallest integer such that $\frac{\lambda}{a+b}\leq n$.

More ambitiously, we may take $\theta_1$, $\theta_2$ as coprime singular inner functions
and consider $\alpha=\theta_1^a \theta_2^b$, $\beta=\theta_1^c \theta_2^d$.
There is then a set of inequalities  that $a,b,c,d$ must satisfy, namely, all lie in $[0,1)$, $(n-1)(a+c) \le 1$, $(n-1)(b+d) \le 1$
with at least one inequality strict,
$n(a+c) \ge 1$, $n(b+d) \ge 1$, and similar inequalities for \eqref{eq:45} and \eqref{eq:46}.
\\

The solution to $Gf=g$ obtained in Proposition \ref{prop:7} below is analogous to the one obtained in \cite{D MSc} for a particular class of symbols $G$ with $\theta=\exp (i\xi)$.

\begin{prop}\label{prop:7}
(i) A solution to $Gf=g$ with $f=(f_1^+,f_2^+)^T \in (H_\infty^+)^2$ and
$g=(g_1^-,g_2^-) \in (H_\infty^-)^2$ is given by
\begin{eqnarray} 
f_1^+ &=& \mu (C_1^{n-1}\overline\alpha^{n-1}-C_1^{n-2}C_2\overline\alpha^{n-2}\beta + \ldots
+ (-1)^{n-1} C_2^{n-1} \beta^{n-1}),\nonumber \\ \label{eq:47}\\
f^+_2 &=& (-1)^n C_2^n \mu \beta^n \overline \theta,\label{eq:48} \\
g_1^- &=& \overline\theta f_1^+,  \\
g_2^- &=& C_1^n \overline\alpha^n \mu,
\end{eqnarray}
where 
\beq\label{eq:51}
\mu=\min\{\theta\overline\beta^{n-1},\alpha\}.
\eeq
(ii) If $\mu=\alpha^n$, then $f=\lambda_1F_+$, where
\beq\label{eq:52}
\lambda_1=\min\{ \alpha^n\beta^n \overline\theta, \alpha \},
\eeq
$F_+=(F_1^+,F_2^+)$ is a corona pair in $(H_\infty^+)^2$ (written $F_+ \in \CP^+$, see\eqref{eq:3.27}), and $g$ is
a corona pair in $(H_\infty^-)^2$ (written $g \in \CP^-$).

(iii) If $\mu=\theta\,\overline\beta^{n-1}$, then $f = \lambda_2 F_+$, where
\beq\label{4.17}
F_+ \in \CP^+\,,\quad\lambda_2 = \min\{ \beta, \theta \overline\alpha^{n-1}\overline\beta^{n-1} \},
\eeq
 and $g \in \CP^-$.
\end{prop}

\beginpf
(i) It is easy to see that $Gf=g$ and $f^+_2 \in H_\infty^+$, $g_2^- \in H_\infty^-$. It remains to prove that
\[
f_{1}^+ \in H_\infty^+, \qquad \overline\theta f_{1}^+ \in H_\infty^-,
\]
for which it is enough to see that $\mu\overline\alpha^{n-1} \in H_\infty^+$
(from \eqref{eq:51} and \eqref{eq:44}) and
$\overline\theta \mu \beta^{n-1} \in H_\infty^-$ (from \eqref{eq:51} and the fact that, if $\mu=\alpha^n$ then $\alpha^n \preceq \theta \overline\beta^{n-1}$).

(ii) If $\mu=\alpha^n$ then
\beq\label{eq:54}
\alpha^n\beta^{n-1} \preceq \theta
\eeq
and
\begin{eqnarray}
f_1^+ &=& \alpha(C_1^{n-1} - C_1^{n-2}C_2 (\alpha\beta) + \ldots + (-1)^{n-1} C_2^{n-1} (\alpha\beta)^{n-1}),
\\
f_2^+ &=& (-1)^n C_2^n \alpha^n \beta^n \overline\theta, \\
g_1^- &=& \overline\theta f_1^+, \\
g_2^- &=& C_1^n.
\end{eqnarray}
Clearly $g=(g_1^-,g_2^-) \in \CP^-$. On the other hand,
\beq
f=\lambda_1 (F_1^+, F_2^+),
\eeq
where $\lambda_1$ is defined by \eqref{eq:52} and $F_1^+, F_2^+ \in H_\infty^+$. If
$\lambda_1 = \alpha^n \beta^n \overline \theta$, then it is clear that $(F_1^+, F_2^+) \in \CP^+$
because $F_2^+ = (-1)^n C_2^n$; if $\lambda_1=\alpha$, then
\begin{eqnarray}\label{eq:60}
F_1^+ &=& C_1^{n-1} - C_1^{n-2}C_2(\alpha\beta)+\ldots + (-1)^{n-1} C_2^{n-1} (\alpha\beta)^{n-1},
\\
\label{eq:61}
F_2^+ &=& (-1)^n (\alpha^{n-1}\beta^n \overline \theta)= (-1)^n C_2^n h^+
\end{eqnarray}
with $h^+ \in H_\infty^+$ because in this case
$\alpha \preceq \alpha^n \beta^n \overline \theta$. We can write
\beq\label{eq:62}
\alpha\beta = (\alpha^{n-1}\beta^n \overline \theta)(\theta \overline\alpha^n \overline \beta^{n-1}) \alpha^ 2 \nonumber\\
= h^+ (\underbrace{\theta \overline\alpha^n \overline\beta^{n-1}}_{\in H_\infty^+ \hbox{ by } \eqref{eq:54}})
\underbrace{\alpha^2}_{\in H_\infty^+},
\eeq
and comparing the expressions \eqref{eq:60} and \eqref{eq:61} for $F_1^+$ and $F_2^+$
respectively, we see that $(F_1^+, F_2^+) \in \CP^+$.

(iii) If $\mu=\theta \overline\beta^{n-1}$, then 
\beq\label{eq:63}
\theta \preceq \alpha^n \beta^{n-1}
\eeq
and
\begin{eqnarray}
f_1^+ &=& \theta(C_1^{n-1} (\overline{\alpha\beta})^{n-1} - C_1^{n-2} C_2  (\overline{\alpha\beta})^{n-2}
+ \ldots + (-1)^{n-1} C_2^{n-1}), \\
f_2^+ &=& (-1)^n C_2^n \beta , \\
g_1^- &=& C_1^n (\overline{\alpha\beta})^{n-1}- C_1^{n-2} C_2 (\overline{\alpha\beta})^{n-2}
+ \ldots + (-1)^{n-1} C_2^{n-1}, \\
g_2^- &=& C_1^n \overline\alpha^n \overline\beta^{n-1} \theta.
\end{eqnarray}
We have $(g_1^-,g_2^-) \in \CP^-$, using \eqref{eq:62} as above. On the other hand,
\[
(f_1^+,f_2^+)= \lambda_2 (F_1^+, F_2^+) \qquad \hbox{with} \quad F_1^+, F_2^+ \in H_\infty^+,
\]
where $\lambda_2$ is defined in \eqref{4.17}. If $\lambda_2 = \beta$, then it is clear that $(F_1^+, F_2^+) \in \CP^+$ since
$F_2^+ = (-1)^n C_2^n$; if
$\lambda_2 = \theta \overline\alpha^{n-1} \overline\beta^{n-1}$, then $\theta\bar\alpha^{n-1} \bar\beta^{n-1} \preceq \beta$, which implies that $\bar\theta\alpha^{n-1}\beta^n\in H_\infty^+$ and
\begin{eqnarray*}
F_1^+ &=& C_1^{n-1} -  C_1^{n-2} C_2 (\alpha\beta) + \ldots + (-1)^{n-1} C_2^{n-1} (\alpha\beta)^{n-1}, \\
F_2^+ &=& (-1)^n C_2^n \, 
\overline\theta \alpha^{n-1}\beta^n.
\end{eqnarray*}
Using the relation $\alpha\beta = (\overline\theta \alpha^{n-1} \beta^n ) ( \theta\overline\alpha^{n-1} \overline\beta^{n-1})\alpha$, where $\theta\overline\alpha^{n-1} \overline\beta^{n-1}\in H_\infty^+$  by \eqref{eq:44}, 
we see as above that $(F_1^+, F_2^+)\in \CP^+$.
\endpf


As a consequence of Theorem \ref{thma:3.12} and Proposition \ref{prop:7}, we have then:

\begin{thm}\label{thm:8}
If \[
G= \begin{bmatrix}
\overline\theta & 0 \\
C_1\overline\alpha + C_2\beta & \theta 
\end{bmatrix},
\]
where $C_1,C_2 \in \CC \setminus\{0\}$ and $\alpha,\beta$ satisfy \eqref{eq:43}--\eqref{eq:46}, then
\[
\ker T_G =\bar\lambda K_\lambda f,
\]
where $f=(f_1^+,f_2^+)$ is defined by \eqref{eq:47}--\eqref{eq:48} and \eqref{eq:51}, and
\[
\lambda = \begin{cases}
\min \{ \alpha^n \beta^n \overline\theta , \alpha\}, & \hbox{if } \alpha^n \beta^{n-1} \preceq \theta, \\
\min \{ \beta, \theta \overline\alpha^{n-1} \overline\beta^{n-1}\}, & \hbox {if } \alpha^n \beta^{n-1} \succeq \theta.
\end{cases}
\]
\end{thm}

The cases $C_1=0$ and $C_2=0$ are rather easier, and we omit them.



\begin{cor}\label{cor4.5}
With the same assumptions as in Theorem \ref{thm:8}
\[
\ker A^\theta_{C_1\bar\alpha+C_2\beta} =\bar \lambda K_\lambda f_1^+.
\]
\end{cor}

\begin{cor}
With the same assumptions as in Theorem \ref{thm:8}, $T_G$ (respectively, $A_h^\theta$) is injective if and only if $\lambda$ is a constant and,
in that case, $T_G$ (respectively, $A_h^\theta$) is invertible.
\end{cor}

\beginpf
The injectivity is a direct consequence of Corollary \ref{cor4.5}. On the other hand, the operator $A_h^\theta$ is equivalent after extension to $T_G$ \cite{BTsk,CP17}; therefore both operators are simultaneously invertible or not. In this case $\det G=1$ and $Gf_+=f_-$ has  a solution $f_\pm \in (H_\infty^\pm)^2$
with $f^\pm\in  \CP^\pm$ and therefore the operator $T_G$ is invertible  \cite{CDR}.
\endpf

\begin{ex}{\rm

Take
$\theta(\xi)=e^{i\xi}$, $\alpha(\xi)=e^{ia\xi}$, $\beta(\xi)=e^{ib\xi}$, $(0<a,b<1)$
and write $
h=C_1 \overline\alpha + C_2 \beta \qquad (C_1,C_2 \in \CC \setminus\{0\})$.
We also write  $\K_\lambda=K_{e_\lambda}$ for $\lambda>0$, where
$e_\lambda(\xi)=e^{i\lambda\xi}$.

Depending on $\alpha$ and $\beta$ there are various possibilities for $\ker A_h^\theta$, some of which we indicate in
Figure \ref{fig:1}, where we have:\\

A: $\K_{a+b-1}e^{i(1-b)\xi}$.

B: $\K_a(C_1-C_2e^{i(a+b)\xi})$.

C: $\K_{1-a-b} (C_1-C_2 e^{i(a+b)\xi})$.

D: $\K_b(C_1 e^{i(1-2b-a)\xi}-C_2e^{i(1-b)\xi})$.

E: $\K_{2a+2b-1}(C_1 e^{i(1-a)\xi}-C_2 e^{i(1-b)\xi})$.


\vbox{ 
\begin{center}
\includegraphics[height=12cm, width=12cm]{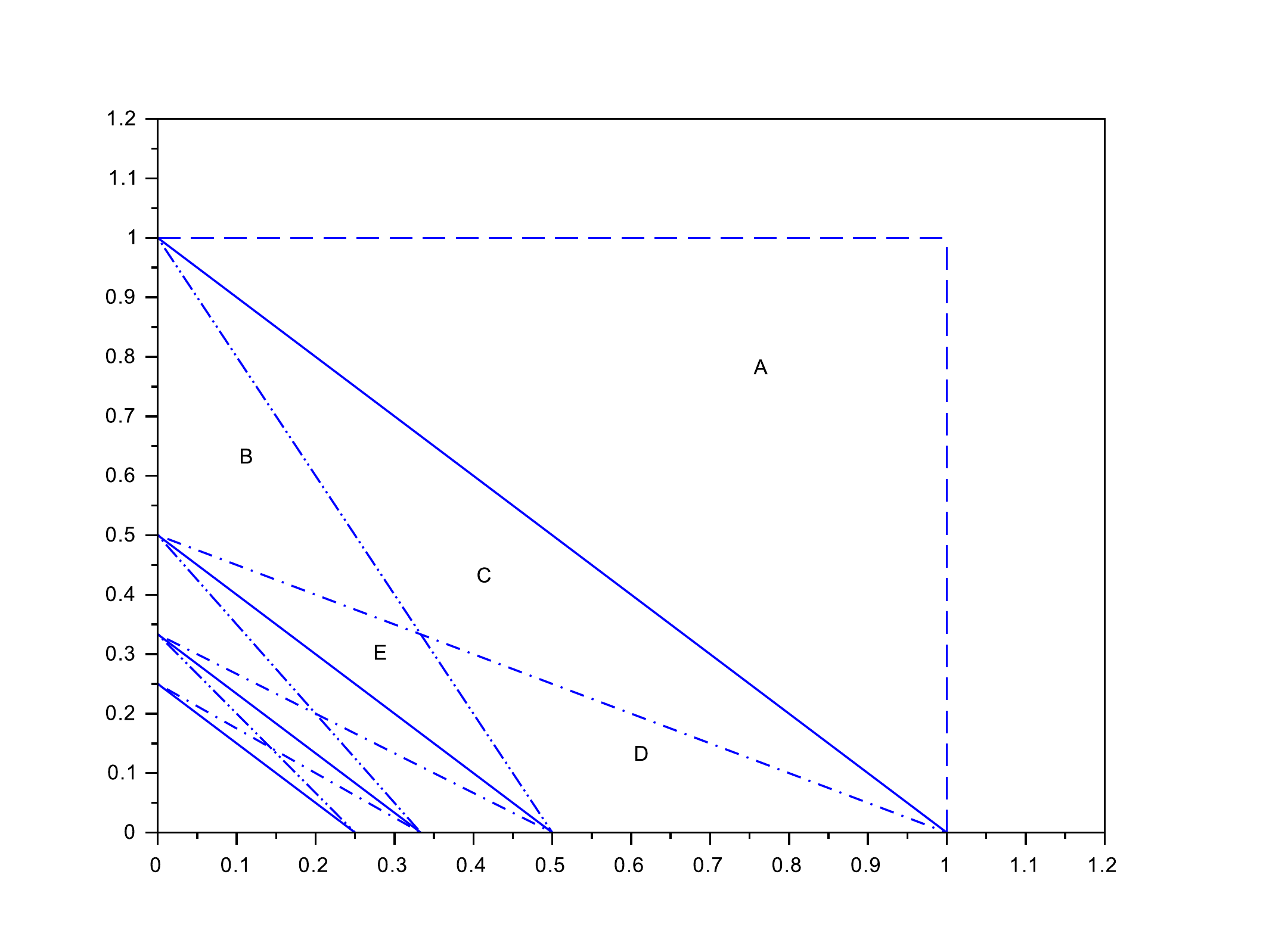}
\captionof{figure}{Dependence of $\ker A_h^\theta$ on $\alpha$ and $\beta$}
\label{fig:1}
\end{center}

} 
}  
\end{ex}

Figure 1 provides a better understanding of the dependence of $\ker A_h^\theta$ on the parameters $a$ and $b$. For instance, one can see that on the lines $a+b=1/n$ the operator is invertible. On the other hand it can easily be verified that the expressions for $\ker A_h^\theta$ ``on the left hand side'' and ``on the right hand side''  of the dotted lines coincide, thus making apparent the continuous dependency of the kernel on the parameters $\alpha$ and $\beta$ across those lines.

\subsection{The second case: $A^\theta_h$ with $h=\overline\alpha R_1^+ + \beta R_2^-$ }

Let
$
h=\overline\alpha R_1^+ + \beta R_2^-
$,
with 
\beq\label{eq:68}
\alpha, \beta \preceq \theta, \qquad \alpha \beta \succ \theta,\qquad
R_1^+ \in \R^+, \qquad R_2^- \in \R^-,
\eeq
where $\R^\pm = \R \cap H_\infty^\pm$. We shall exclude the degenerate cases $R_1^+=0$ and $R_2^-=0$.

If $\alpha=\beta=\theta$, then the operators $A^\theta_h$ are the so-called
finite-rank truncated Toeplitz operators of Type I \cite{GuKang,Bess}.

We assume here that $\theta$ is an inner function that is not a finite Blaschke product (written $\theta \not\in \FBP$),
otherwise the matrix $G$ would be rational, and also that $\alpha,\beta, \alpha\beta/\theta \not\in \FBP$. 
Note that if $\alpha,\beta \in \FBP$, then $h \in \R$; this case was studied in \cite{CP16}.

We have:

\beq\label{eq:70}
G= \begin{bmatrix}
\overline\theta & 0 \\ \overline\alpha R_1^+ + \beta R_2^- & \theta
\end{bmatrix},
\eeq
and $Gf=g$ holds with
\beq\label{eq:71}
f = \begin{bmatrix}
\alpha \\ -\frac{\alpha\beta}{\theta}R_2^- 
\end{bmatrix}, \qquad
g= \begin{bmatrix}
\overline\theta\alpha \\ R_1^+  
\end{bmatrix}.
\eeq
Let 
\beq
R_1^+ = \frac{N_1}{D_1^+},
\eeq
where $N_1$ and $D_1^+$ are polynomials without common zeroes, with
$\deg N_1 \le \deg D_1^+ = n_1$, such that all zeroes of $D_1^+$ are in $\CC^-$, and
\beq
F_2^-=\frac{N_2}{D_2^-},
\eeq
where $N_2$ and $D_2^-$ are polynomials without common zeroes with
$\deg N_2 \le \deg D_2^-=n_2$, such that all zeroes of $D_2^-$ are in $\CC^+$.
Condition \eqref{eqa:3.16}
is satisfied in this case, by Lemma \ref{lem:14} below. In fact we have
\[
\begin{array}{rl}
&\left[\alpha \quad -\frac{\alpha\beta}{\theta}R_2^-  \right]
\begin{bmatrix}
\phi_{1+} \\ \phi_{2+}
\end{bmatrix}
=
[\overline\theta \alpha \quad R_1^+]
\begin{bmatrix}
\phi_{1-} \\ \phi_{2-}
\end{bmatrix}
\\
\iff & \alpha \phi_{1+}-\frac{\alpha\beta}{\theta}R_2^-  \phi_{2+} = \overline\theta \alpha\phi_{1-} + R_1^+\phi_{2-}
\\
\iff & \underbrace{\frac{\alpha\beta}{\theta}}_{\not\in \FBP}
\underbrace{D_1^+ \overline{D_2^-}}_{p_1}
\underbrace{\left( \frac{\theta}{\beta} \frac{D_2^-}{\overline{D_2^-}}\phi_{1+}-\frac{N_2}{\overline{D_2^-}}\phi_{2+}\right)}_{\in H_2^+}
=
\underbrace{D_2^- {\overline{D_1^+}}}_{p_2}
\underbrace{\left( \frac{\overline\theta \alpha D_1^+}{\overline{D_1^+}} \phi_{1-}+\frac{N_1}{\overline{D_1^+}}\phi_{2-}\right)}_{\in H_2^-},
\end{array}
\]
so both sides of the equation must be equal to zero.

We have the following left inverses for $f$ and $g$:
\beq
\tilde f= (\overline\alpha, 0), \qquad \tilde g =(\theta\overline\alpha, \quad 0),
\eeq
so $\HH_\pm$ are defined by
\begin{eqnarray}
\HH_+ &=& \{(\psi_{1+},\psi_{2+})  \in (H_2^+)^2: R_2^- \beta \overline\theta \psi_{1+} \in H_2^+ \},\nonumber\\
\HH_- &=&  \{(\psi_{1-},\psi_{2-})  \in (H_2^-)^2: \theta\overline\alpha R_1^+ \psi_{1-} \in H_2^- \}.
\label{eq:73}
\end{eqnarray}
From \eqref{eq:73}, we have that, for $\phi_+ \in H_2^+$:
\beq\label{eq:76}
R_2^- \overline\theta \beta \psi_{1+} = \phi_+
\iff
\frac{N_2}{D_2^-} \psi_{1+} = \frac{\theta}{\beta} \phi_+ 
 \iff 
\frac{N_2}{\overline{D_2^-}} \psi_{1+} = \frac{\theta}{\beta} \frac{D_2^-}{\overline{D_2^-}} \phi_+.
\eeq
Now let $\ds \left( \frac{N_2}{\overline {D_2^-}}\right)_i$ denote the inner factor in an inner--outer
factorization of $\ds \frac{N_2}{\overline {D_2^-}}$. Since
$\ds \frac{\theta}{\beta} \frac{D_2^-}{\overline{D_2^-}}$  is inner and there are no common zeroes for $N_2$ and $D_2^-$ we have
\beq\label{eq:77}
\gamma_2: = \gcd \left\{ \left( \frac{N_2}{\overline {D_2^-}}\right)_i, \frac{\theta}{\beta} \frac{D_2^-}{\overline{D_2^-}}  \right\}
=
\gcd\left \{  \left( \frac{N_2}{\overline {D_2^-}}\right)_i, \frac{\theta}{\beta} \right \},
\eeq
and it follows from \eqref{eq:73} and
\eqref{eq:76} that $( \psi_{1+},\psi_{2+}) \in \HH_+$ if and only if
$\psi_{1+},\psi_{2+} \in H_2^+$ and
\beq
\psi_{1+} \in \frac{\theta}{\beta \gamma_2} \frac{D_2^-}{\overline{D_2^-}} H_2^+.
\eeq
Analogously,  defining
\beq\label{eq:80}
\gamma_1 = \gcd \left \{ \left( \frac{\overline{N_1}}{D_1^+} \right)_i, \frac{\theta}{\alpha} \frac{\overline{D_1^+}}{D_1^+} \right \}
=
\gcd \left\{ \left( \frac{\overline{N_1}}{D_1^+} \right)_i, \frac{\theta}{\alpha}  \right \},
\eeq
we have that 
$( \psi_{1-},\psi_{2-}) \in \HH_-$ if and only if
\beq
\psi_{1-},\psi_{2-} \in H_2^-\quad {\rm  and} \quad\psi_{1-} \in \frac{\alpha}{\theta} \gamma_1  \frac{D_1^+}{\overline{D_1^+}} H_2^-.
\eeq
Therefore, from \eqref{eqa:3.20}, $\K$ is defined by the equation
\beq
[\overline\alpha \quad 0] \begin{bmatrix}
\frac{\theta}{\beta\gamma_2} \frac{D_2^-}{\overline{D_2^-}}\phi_+ \\ \psi_{2+}
\end{bmatrix}
= 
[\theta \overline\alpha \quad 0]
\begin{bmatrix}
\frac{\alpha\gamma_1}{\theta}  \frac{D_1^+}{\overline{D_1^+}}\phi_- \\ \psi_{2-}
\end{bmatrix}
\eeq
with $\phi_\pm$, $\psi_{2\pm} \in H_2^\pm$, i.e.,
\beq
\K = \left( \frac{\theta}{\alpha\beta} \overline{\gamma_2} \frac{D_2^-}{\overline{D_2^-}} H_2^+ \right) \cap
\left( \gamma_1 \frac{D_1^+}{\overline{D_1^+}} H_2^- \right).
\eeq
Consequently,
\begin{eqnarray*}
\ker T_G &=& \K f = \K \begin{bmatrix}
\alpha \\ - \frac{\alpha\beta}{\theta}R_2^-
\end{bmatrix}
=
K\frac{\alpha\beta}{\theta}\ \begin{bmatrix}
\theta/\beta \\ -R_2^- 
\end{bmatrix}
\\
&=& \underbrace{H_2^+ \cap
\left( \gamma_1\gamma_2 \frac{D_{1^+}}{\overline{D_1^+}}\frac{\overline{D_2^-}}{D_2^-}\frac{\alpha\beta}{\theta} H_2^- \right)}_{\ker T_{\eta^{-1}}}
\underbrace{\begin{bmatrix}
\overline{\gamma_2}\frac{\theta}{\beta} \frac{D_2^-}{\overline{D_2^-}} \\
- \overline{\gamma_2} \frac{N_2}{\overline{D_2^-}}
\end{bmatrix}
}_{\in (H_\infty^2)^2 \hbox{ by } \eqref{eq:77} \hbox{ and } \eqref{eq:80}}
\end{eqnarray*}
with
\beq\label{4.47}
\eta = \gamma_1 \gamma_2 \frac{D_{1^+}}{\overline{D_1^+}}\frac{\overline{D_2^-}}{D_2^-} \frac{\alpha\beta}{\theta}.
\eeq

We have thus shown:

\begin{thm}
If $g=\overline\alpha R_1^+ + \beta R_2^-$, where $\alpha,\beta$ are inner
functions and \eqref{eq:68} is satisfied, then, for $G$ defined by \eqref{eq:70},
we have
\beq
\ker T_G = \ker T_{\eta^{-1}} \begin{bmatrix}
\overline{\gamma_2}\frac{\theta}{\beta} \frac{D_2^-}{\overline{D_2^-}} \\
- \overline{\gamma_2} \frac{N_2}{\overline{D_2^-}}
\end{bmatrix}\quad {\rm and} \quad \ker A_g^\theta = 
\overline{\gamma_2}\frac{\theta}{\beta} \frac{D_2^-}{\overline{D_2^-}}\,\ker T_{\eta^{-1}}\,,
\eeq 
with the notation above and $\eta$ given by \eqref{4.47}.
\end{thm}

\noindent{\bf Remark}.
Note that $ 
\overline{\gamma_2}\frac{\theta}{\beta} \frac{D_2^-}{\overline{D_2^-}}$ is an inner function. Also note that 
if 
\[
\left(\overline\gamma_1 \frac{\theta}{\alpha}\frac{\overline{D_1^+}}{D_1^+}\right)
\left(\overline\gamma_2  \frac{\theta}{\beta}\frac{D_2^-}{\overline{D_2^-}}\right)\preceq \theta
\]
 then $\ker T_{\eta^{-1}}$ is a model space, and in that case we have $\ker T_G=\{0\}$ (and $\ker A^\theta_g=\{0\})$
if and only if $\eta$ is a constant.

\subsection*{Some auxiliary results}

\begin{lem}\label{lem:12}
Let $\alpha$ be an inner function, $\alpha \not\in \FBP$, and let $p_1,p_2$ be polynomials. Then
$\alpha p_1 H_2^+ \cap p_2 H_2^- = \{0\}$.
\end{lem}
\beginpf
If there are $\phi_\pm \in H_2^\pm$ such that
$
\alpha p_1 \phi_+ = p_2 \phi_-
$,
then both sides of this equation must be equal to a polynomial $p$, by an easy generalization of Liouville's Theorem,
and we have $
\alpha \phi_+ = \frac{p}{p_1} \in H_2^+, \, \hbox{i.e.,} \,\, \overline\alpha \frac{p}{p_1} = \phi_+ \in H_2^+$.
Therefore, $\ds \overline\alpha p/p_1 = P^+ \left(\overline\alpha p/p_1\right)$ must be rational, which
is impossible since $\alpha \not\in \FBP$.
\endpf

\begin{cor}
$\alpha \HH_2^+ \cap \HH_2^- = \{0\}$,
where $\HH_2^\pm = (\xi \pm i)H_2^\pm$.
\end{cor}

As a consequence of Lemma \ref{lem:12} we can obtain the following generalization:

\begin{lem}\label{lem:14}
Let $\alpha$ be an inner function, $\alpha\not\in \FBP$, and let
$\ds R_1^+ = \frac{N_1}{D_1^+}$, $\ds R_2^- = \frac{N_2}{D_2^-}$,
where 
$N_1$, $D_1^+$, $N_2$, $D_2^-$ are polynomials such that
$N_1$ and $D_1^+$ (and similarly $N_2$ and $D_2^-$) have no common zeroes,
all zeroes of $D_1^+$ (respectively,  $D_2^-$) are in $\CC^-$ (respectively, $\CC^+$)
and $\deg N_1 \le \deg D_1^+ = n_1$ and $\deg N_2 \le \deg D_2^- = n_2$.
Then 
\beq\label{eq:85}
\alpha R_2^- H_2^+ \cap R_1^+ H_2^- = \{0\}.
\eeq
\end{lem}
\beginpf
If $\phi_\pm \in H_2^\pm$ and
\[
\alpha \frac{N_2}{D_2^-} \phi_+ = \frac{N_1}{D_1^+}\phi_-,
\]
then
\[
\alpha
\underbrace{D_1^+ \overline{D_2^-}}_{p_1}
 \underbrace{\left( \frac{N_2}{\overline{D_2^-}} \phi_+\right)}_{\in H_2^+} 
=
\underbrace{D_2^- \overline{D_1^+}}_{p_2}
 \underbrace{\left(\frac{N_1}{\overline{D_1^+}}\phi_-\right)}_{\in H_2^-}  = 0
\]
by Lemma \ref{lem:12}, and \eqref{eq:85} follows.
\endpf

\subsection*{Acknowledgements}

This work was partially supported by FCT/Portugal through\\  UID/MAT/04459/2019.

\end{document}